\documentclass[reqno,12pt]{amsart}
\usepackage{amsmath,amssymb,amsthm}
\usepackage{float}
\usepackage{geometry}
\usepackage{tabularx}
\usepackage[dvipsnames]{xcolor}
\usepackage[
         breaklinks=true,
		colorlinks=true,
		linkcolor=blue,
		citecolor=ForestGreen,
		urlcolor=cyan
		]{hyperref}
\usepackage{mathtools}\usepackage{soul}

\usepackage{subfigure}\usepackage{tikz}

\newtheorem{lemma}{Lemma}[section]

\newtheorem{proposition}{Proposition}[section]

\newtheorem{theorem}{Theorem}[section]

\theoremstyle{definition}
\newtheorem{definition}{Definition}[section]

\theoremstyle{remark}

\theoremstyle{remark}
\newtheorem{remark}{Remark}[section]
\numberwithin{equation}{section}

\newcommand{\C}{{\mathbb C}}

\newcommand{\R}{{\mathbb R}}

\newcommand{\la}{{\mathcal L}_a}
\definecolor{blu}{rgb}{0,0,1}

\newcommand{\supp}{\rm{supp}}

\newcommand{\rd}{{\mathbb R}^d}

\allowdisplaybreaks
\newtheorem{itheorem}{{\bf Theorem}}[section]

\title[The heat equation with inverse square potential]{On the hardy-H\'enon heat equation with an inverse square potential}

\author[D.  G.  Bhimani]{Divyang G. Bhimani}
\address{Department of Mathematics, Indian Institute of Science Education and Research Pune, Dr. Homi Bhabha Road, Pune 411008, India}
\email{divyang.bhimani@iiserpune.ac.in}

\author[S. Haque] {Saikatul Haque}
\address{Department of Mathematics, Indian Institute of Science Education and Research Tirupati, Andhra Pradesh 517619, India}
\email{saikatul@iisertirupati.ac.in}

\author[M. Ikeda]{Masahiro Ikeda}
\address{Faculty of Science and Technology,
Keio University, 
3-14-1 Hiyoshi, Kohoku-ku, Yokohama, 223-8522, Japan \& Center for Advanced Intelligence Project
RIKEN, Japan.}
\email{masahiro.ikeda@keio.jp, masahiro.ikeda@riken.jp}
\begin{document}
\subjclass[2000]{35K05,  35K67, 35B30, 35K57}
\keywords{ Hardy-H\'enon equation,  Inverse square potential,  Dissipative estimate,  well-posesness,   finite time blow up}

\begin{abstract} We study Cauchy problem for the Hardy-H\'enon parabolic equation with an inverse square potential, namely, 
\[\partial_tu  -\Delta u+a|x|^{-2} u= |x|^{\gamma} F_{\alpha}(u),\]
where $a\ge-(\frac{d-2}{2})^2,$ $\gamma\in \R$,  $\alpha>1$  and $F_{\alpha}(u)=\mu |u|^{\alpha-1}u, \mu|u|^\alpha$ or $\mu u^\alpha$, $\mu\in \{-1,0,1\}$.  We  establish sharp fixed time-time decay estimates   for heat semigroups $e^{-t (-\Delta + a|x|^{-2})}$  in  weighted Lebesgue spaces. This may be of independent interest.   

As an application,   we establish   local well-posedness  in  scale subcritical  and critical   weighted Lebesgue spaces and small data global existence in critical weighted Lebesgue spaces. Further, under certain conditions on $\gamma$ and $\alpha,$ we show that  local solution  cannot be extended to global  one for certain initial data in the subcritical regime.   Thus, finite time blow-up in  the subcritical Lebesgue space norm is exhibited.  We also demonstrate  nonexistence of local positive weak solution (and hence failure of local well-posedness) in supercritical case for $\alpha>1+\frac{2+\gamma}{d}$ the Fujita exponent. 
\end{abstract}
\maketitle
\section{Introduction}
\subsection{Fixed-time  estimates for heat semigroup  $e^{-t (-\Delta + a |x|^{-2})}$}
Consider the linear heat equation associated with the inverse square potential,  namely
\begin{equation}\label{dsr}
\begin{cases}  \partial_t u(t,x) + \mathcal{L}_au(t,x)=0\\
u(0,x)= u_0(x)
\end{cases} (t,x) \in \mathbb R^+ \times \mathbb R^d, 
\end{equation}
where $u(t,x)\in \mathbb C.$   In this paper,  we  assume that $a\ge a_*:=-(\frac{d-2}{2})^2, d\geq 2,$ unless it is explicitly specified.   The
 Schr\"odinger operator with inverse square potentials$$\la=-\Delta+ a |x|^{-2}$$
 is initially defined with domain $C^\infty(\rd\setminus\{0\})$. 
Then it is extended as an unbounded operator in weighted Lebesgue space $L_s^q(\rd)$ 
that generates a positive semigroup $\{e^{-t\la}\}_{t\ge0}$ provided $1<q<\infty$ and $\sigma_-<\frac{d}{q}+s<\sigma_++2,$ where $\sigma_-$, $\sigma_+$  defined  by
\begin{equation}\label{sigma}
\sigma_\mp=\sigma_\mp(d,a):=\frac{d-2}{2}\mp\frac{1}{2}\sqrt{(d-2)^2+4a}
\end{equation} are the roots of $s^2 - (d - 2)s - a = 0$, see \cite[Theorems 3.2, 3.3]{metafune2021maximal}. Here,   the weighted Lebesgue space $L_s^q(\rd)$   is defined  by the norm $\| f\|_{L^q_s}: = \| |\cdot|^{s}f\|_{L^q} \  (s\in \mathbb R).$

The study of  $\mathcal{L}_a$ is motivated from physics and mathematics spanning areas such as combustion theory, the Dirac equation with Coulomb potential, quantum mechanics and the study of perturbations of classic space-time metrics.  See e.g.  \cite{JLVazquez, Kalf,  Burq} and the references therein.

The aim of this article  is to understand the dynamics of solutions of  Hardy-H\'enon heat equations  \eqref{dsr} and \eqref{heat equation} when a  singular potential is present,  in light  of the research programme initiated by  Zhang \cite{Zhang2001},  Pinsky  \cite{Pinsky2009, Pinsky1997},  Ioku et al.  in \cite{ioku2016lplq, Ioku2019},   Ishige \cite{Ishige2007} and Ishige-Kawakami in  \cite{Ishige2020},   and Bhimani-Haque  \cite{Divyang_Saikatul_arxiv} (cf.  \cite{bhimani2021,  bhimani2023, bhimani2023heat}).  We also note that there is  a extensive literature on Hardy-H\'enon heat equation  without potential,  i.e.  \eqref{heat equation}  with $a=0,$  we refer to recent work of   Chikami et al.   in \cite{chiiketanitay_arXiv,  chikami2022optimal} and the references therein,  see also Remark \ref{HR}.

     We begin by stating our dissipative estimates in weighted Lebesgue spaces in the following theorem.

\begin{theorem}\label{wle} Let $\sigma_-,\sigma_+$ be as defined in \eqref{sigma}. Let $s_1,s_2\in \R$ and $q_1,q_2\in (1,\infty)$. Then 
\begin{equation}\label{main est}
   \|e^{-t\mathcal{L}_a}f\|_{L^{q_2}_{s_2}}\le Ct^{-\frac{d}{2}\left(\frac{1}{q_1}-\frac{1}{q_2}\right)-\frac{s_1-s_2}{2}}\|f\|_{L^{q_1}_{s_1}}\qquad\forall \ t>0,\ \forall \ f\in L^{q_1}_{s_1}(\R^d)
\end{equation}
if and only if
\begin{equation}\label{Condiexp}
\sigma_-<\frac{d}{q_2}+s_2\le \frac{d}{q_1}+s_1<\sigma_++2,
\end{equation}and
\begin{equation}\label{reg}
s_2\le s_1.
\end{equation}
\end{theorem}

\begin{remark}\label{HR} Theorem \ref{wle} deserve several comments.
\begin{enumerate}
\item \label{HR1} The case $a=0$:   In this case $e^{-t\mathcal{L}_0}f=e^{t\Delta}f=k_t\ast f$ (where $k_t:=t^{-d/2}\exp(-\frac{|\cdot|^2}{4t})$) and $\sigma_-=0,\sigma_++2=d$.
\begin{itemize}
\item Subcase $s_1,s_2=0$: The sufficiency part \eqref{main est} is a consequence of Young's convolution inequality.   See \cite[Lemma 3.1]{Miao08}. This argument  holds even if we replace strict inequities in \eqref{Condiexp} by  equalities and thus $q_1,q_2$ can take the extreme values $1,\infty$. 
\item Subcase $s_1$ or $s_2\in\R\setminus\{0\}$: For $q_1\le q_2$, this  is due to Chikami-Ikeda-Taniguchi  \cite[Lemma 2.1]{chikami2022optimal}.    Theorem \ref{wle}
 removes the assumption $q_1\le q_2$ in \cite[Lemma 2.1]{chikami2022optimal}. See also \cite[Section 3]{tayachi2023new}  and \cite{slimene2017well,tayachi2020uniqueness,tsutsui2011navier}.
 \item Recently  Chikami-Ikeda-Taniguchi-Tayachi established similar estimate in Lorentz spaces. Their result does not cover the case $\frac{d}{q_2}+s_2= \frac{d}{q_1}+s_1$, $s_2<s_1$, see  \cite[Proposition 3.1]{chiiketanitay_arXiv}.
\end{itemize}
\item The case $a\in [a_*,\infty)$:
\begin{itemize}
\item Subcase $s_1,s_2=0$: In this subcase,
the sufficiency part \eqref{main est}
 is due to Ioku-Metafune-Sobajima-Spina \cite[Theorem 5.1]{ioku2016lplq}.     However,  their method of  proof  is different than ours,  which rely  on  embedding theorems and interpolation techniques. The
 
\item Subcase $s_1 $ or $ s_2\in\R \setminus \{ 0 \} $:  In this case,  both necessity and sufficiency part of Theorem \ref{wle} is new.   This is the main contribution of this article.
\end{itemize}

\item The power of $t$ in \eqref{main est} is optimal which follow  by a standard scaling argument, see Lemma \ref{lem:scaling}.


\item\label{symm} Using Symmetry (in $x,y$ variable) of heat kernel $g_a(t,x,y)$ (see Subsection \ref{HK}) associated with the operator $e^{-t\la}$, it follows by duality and the relation $\sigma_++2=d-\sigma_-$ 
that \eqref{main est} holds for $(q_1,s_1,q_2,s_2)$ if and only if  \eqref{main est} holds for $(q_2',-s_2,q_1',-s_1)$ (here $q_j'$ is  the H\"older conjugate of $q_j$).

\item For  $s_1=-\sigma_-$,   Theorem \ref{wle} holds even for end point cases $q_1\in\{1,\infty\}$ (hence allowing equality in the last strict inequality in \eqref{Condiexp}).    For  $s_2=\sigma_-$,   Theorem \ref{wle} holds even for end point cases $q_2\in\{1,\infty\}$ (hence allowing equality in the first strict inequality in \eqref{Condiexp}). 


\item \label{HRD} It is indispensable  to  consider  weighted Lebesgue spaces  in Theorem \ref{wle}   in 
 order to  treat H\'enon potential $|x|^{\gamma} \ (\gamma>0)$ while establishing well-posedness for \eqref{heat equation}.
\end{enumerate}
\end{remark}

\begin{remark}  We briefly comment on the novelty and key ideas  of the  proof of Theorem \ref{wle}.
\begin{enumerate}
    \item The classical  heat semigroup  $e^{t\Delta}$ (case $a=0$) is a  convolution operator with  the  Gaussian kernel. Hence, for non-weighted case $s_1=s_2=0$, the analogue of Theorem \ref{wle} (sufficiency part) is a consequence of Young's inequality and scaling argument  (see Remark \ref{HR} \eqref{HR1}). This does not work if $s_1$ or $s_2$ is non-zero, and we treat this case along with the case $a\neq0$.
    \item The operator   $e^{- \mathcal{L}_a}$ for the case $a\neq 0$ is not a convolution operator\footnote{We only need to focus on estimating in the case $t=1$, see Lemma \ref{lem:scaling}}.  This raises some difficulties,  and  method used for the classical case is not applicable.  In fact, in this case,  we rely on  the upper-bound of heat kernel  $g_a(1,x,y)$ (Theorem \ref{heat kernel}), which contains the singularity at the origin (though it has rescaled Gaussian multiplier present). To overcome this challenge, we split the area of integration in $x$-$y$ plane into four regions (depending on whether $x,y$ are near  or away from origin); and  use different techniques for each region  to achieve the desired estimate. 
    \item  We   use Holder's and Young's inequalities with  cleverly   chosen Lebesgue spaces when the middle inequality in \eqref{Condiexp} is strict. On the other hand, when the middle inequality in \eqref{Condiexp} is equality, Lebesgue spaces alone are not enough.  In  this  case, we  use  suitable Lorentz spaces  to complete the argument.
\end{enumerate}  
\end{remark}

\subsection{Hardy-H\'enon equations (HHE) with inverse-square potential}
We  consider \eqref{dsr}  with an inhomogeneous power  type nonlinearity:
\begin{equation}\label{heat equation}
\begin{cases}
\partial_tu(t,x)+ \mathcal{L}_au(t,x)=|x|^{\gamma}F_{\alpha}(u(t,x))\\
u(x, 0)= u_0(x)
\end{cases}
 (t,x) \in [0,T) \times \rd,    
\end{equation}
where $\gamma\in \R$, 
$T\in (0,\infty],$ and $\alpha>1$ and $u(x,t)\in \mathbb R$ or $u(x,t)\in \mathbb C.$  We assume that  the non-linearity  function $F_{\alpha}:\C\rightarrow \C$  satisfies the following conditions:
\begin{eqnarray}\label{Fcon}
\begin{cases}  |F_{\alpha}(z)-F_{\alpha}(w)|\le C_0(|z|^{\alpha-1}+|w|^{\alpha-1})|z-w| & for \ z, w \in \mathbb C\\
F_\alpha(0)=0.
\end{cases}
\end{eqnarray}
The typical examples of $F_\alpha$ would be \[ F_{\alpha}(z)=\mu |z|^{\alpha-1}z,  \ \mu |z|^\alpha \text{ or }  \   \mu z^\alpha \quad (\mu\in \mathbb R).
\] 
 The  potential  $|x|^{\gamma}$  is called H\'enon type if $\gamma>0$ and is called {Sobolev} type if $\gamma<0$.  The equation \eqref{heat equation} with $\gamma<0$ is known as a 
{\sl Hardy parabolic equation}, while that with $\gamma>0$ 
is known as a {\sl H\'enon parabolic equation}.  Equation (\ref{heat equation}) is called Hardy-H\'enon parabolic equation with an inverse square potential. 
The elliptic part of \eqref{heat equation} when $a=0$,  i.e.
\[-\Delta u= |x|^{\gamma}|u|^{\alpha-1}u \]
was proposed by H\'enon \cite{henon1973numerical}  as a model to study the rotating stellar systems  and has been extensively  studied in scientific community,  see e.g. \cite{Nassif}.

The equation (\ref{heat equation}) is invariant under the following scale transformation:
\[
u_\lambda(t,x) := \lambda^\frac{2+\gamma}{\alpha-1}u(\lambda^2 t, \lambda x),\quad \lambda>0.
\]
More precisely, if $u$ is a solution to \eqref{heat equation}, then so is $u_\lambda$
with the rescaled initial data $\lambda^\frac{2+\gamma}{\alpha-1}u_0(\lambda x)$.  
Then the following identity holds 
\[
\|u_\lambda(0)\|_{L^{q}_s} = \lambda^{-s+ \frac{2+\gamma}{\alpha-1}- \frac{d}{q}}
\|u_0\|_{L^{q}_s},
\quad \lambda>0.
\]
Hence, if $q$ and $s$ satisfy 
\begin{equation*}
s + \frac{d}{q} = \frac{2+\gamma}{\alpha-1},
\end{equation*}
then the identity $
\|u_\lambda(0)\|_{L_s^q} 
=
\|u_0\|_{L_s^q}$ holds for any $\lambda>0$, 
i.e., the norm $\|u_\lambda(0)\|_{L^{q}_{s}} $ is invariant with respect to $\lambda$. 
Denote
\begin{equation}\label{T}
\tau= \tau (q,s,d):=s + \frac{d}{q} \quad\text{and}\quad\tau_c=\tau_c(\gamma,\alpha):=\frac{2+\gamma}{\alpha-1}.
\end{equation} 
 We say   Cauchy problem \eqref{heat equation} scale 
\begin{eqnarray}\label{cri}
L_s^q-\begin{cases}
subcritical &  \text{if}  \  \   \tau<\tau_c\\
 critical &  \text{if}  \  \   \tau=\tau_c\\
supercritical & \text{if}  \  \  \tau>\tau_c
\end{cases}.
\end{eqnarray}
\begin{remark}\label{cr0} For $\tau=\tau_c,$ we get $s=\frac{2+\gamma}{\alpha-1}- \frac{d}{q}=:s_c(q,\gamma,\alpha,d)$ (often denoted by $s_c$  for shorthand). 
 In particular,  when $s=s_c=0,$ $\gamma\ge-2$, we have $q=q_c:=\frac{d(\alpha-1)}{2+\gamma}=\frac{d}{\tau_c}.$  So $L^{\frac{d(\alpha-1)}{2+\gamma}}(\rd)$ is the critical Lebesgue space without weight.
\end{remark}
We recall  the notion of  well-posedness  in the sense of Hadamard.
\begin{definition}\label{well}[well-posedness]
Let $T\in(0,\infty], s\in \mathbb R$, $1\leq q \leq \infty$ and $\sigma_-<\tau<\sigma_++2$.
\begin{itemize}
\item[--]  We say that $u$ is an $L_s^q$-integral solution on $[0,T)$ to (\ref{heat equation}) if $u\in C([0,T); L^q_s(\R^d))$ and satisfies
\begin{equation}\label{integraleq}
   u(t)=e^{-t\mathcal{L}_a}u_0+\int_0^te^{-(t-\tau)\mathcal{L}_a}[|\cdot|^{\gamma}F_\alpha(u(\tau))]d\tau
\end{equation}
for any $t\in [0,T)$. Maximum of such $T$ is denoted by $T_m$.
\item[--]
Let $X,Y \subset \mathcal{S}'(\rd)$ be Banach spaces. Then
 \eqref{heat equation} is called locally
 well-posed (in short LWP) from  $X$ to $Y$ if, for each bounded  $B\subset X$, there exist $T>0$ and a Banach space  $X_{T} \hookrightarrow C([0,T], Y)$ so that
\begin{enumerate}
\item[(a)] for all $u_0\in B$,  \eqref{heat equation}  has a unique integral solution $u\in X_T$ 
\item[(b)]  $u_0\mapsto u$ is continuous from $(B, \|\cdot\|_{X})$ to $C([0,T],Y).$
\end{enumerate}
If $X=Y$ we say  \eqref{heat equation} is locally well-posed {in  $X$}.
 If $T=\infty$, 
 then we say \eqref{heat equation} is {globally} well-posed in $X$.
\end{itemize}
\end{definition}

\begin{remark} \label{hr}We briefly mention some history  on several facets of  \eqref{heat equation}.  We define   Fujita exponent by
\[\alpha_F=\alpha_F(d, \gamma,a)= 
1+\frac{(2+\gamma)^+}{\sigma_++2}, \qquad (\text{here }(2+\gamma)^+=\max(2+\gamma,0))\]
which is often known to divide  the existence and nonexistence of positive global solutions.
\begin{enumerate}\item By taking  $a=\gamma=0$  and $F_{\alpha}(z)=z^{\alpha}$ in \eqref{heat equation},  we get classical heat equation
\begin{eqnarray}\label{d1}
\partial_tu - \Delta u= u^{\alpha}, \quad u(0)=u_0.
\end{eqnarray}
We recall following  known results for \eqref{d1}: 
\begin{enumerate}
\item \label{hr1} Let $q_c$ be as in  Remark \ref{cr0}. If $q\geq q_c$ and $q>1$ or $q>q_c$ and $q\geq 1,$ Weissler \cite{weissler1980local}  proved  the existence of a unique local solution  $u \in C([0, T), L^q(\rd)) \cap L^{\infty}_{loc} (0, T],  L^{\infty}(\rd)).$ Later on,  Brezis-Cazenave \cite{brezis1996nonlinear} proved the unconditional uniqueness of Weissler's solutions.
\item\label{hr2}  If $q<q_c$, there are indications that there exists no (local) solution in any reasonable weak sense, see \cite{weissler1980local, brezis1996nonlinear, weissler81}. Moreover, it is known that uniqueness
is lost for the initial data  $u_0=0$ and for  $1+ \frac{1}{d}< q < \frac{d+2}{d-2}$, see \cite{haraux1982non}.
\item\label{hr4} Fujita \cite{Fujita} proved, for  $1< \alpha <\alpha_F(d, 0,0),$ \eqref{d1} has no global solution (i.e. every solution blows up in finite time in  $L^{\infty}-$norm), whereas  for $\alpha> \alpha_{F}(d, 0,0),$  classical solution  is global for small data.
\end{enumerate}
\item Taking $a=0,  F_{\alpha}(z)=z|z|^{\alpha-1}$ in  \eqref{heat equation},  we  get classical Hardy-H\'enon heat equation
\begin{equation}\label{d2}
\partial_tu- \Delta u= |x|^\gamma  |u|^{\alpha-1}u, \quad u(0)=u_0.
\end{equation}
 In this case,  Chikami et al.   in \cite{chikami2022optimal} introduced weighted Lebesuge space $L^q_s(\rd)$ to treat potential $|x|^{\gamma},$ and establish well-posedness results.  Later,  Chikami et al.  in  \cite{chiiketanitay_arXiv} generalize these results in weighted Lorentz spaces.    In this paper,  we could  establish analogue of these results in the presence of potential,  i.e.  for \eqref{heat equation} with $a\neq 0$ and relaxed conditions on other parameters $\gamma,\alpha,q,s$.  See Remarks \ref{wgr} and \ref{npsr} below.

\item Several authors considered \eqref{heat equation} with some mild restriction on external potential:
\begin{eqnarray}\label{d3}
\partial_t u  -\Delta u -V (x)u= b(x)u^{\alpha},  \quad u(0)=u_0, 
\end{eqnarray}
 and showed sharp contrast between existence of classical global solution and finite time blow-up in $L^\infty-$norm by finding appropriate Fujita exponent.  We recall some of them here:
\begin{enumerate}
\item \label{p97}  Let $V(x)=\frac{a}{|x|^2}$ and $b \in C^{\beta}(\R)$ ($\beta\in(0,1]$) with $b(x)\sim|x|^\gamma$ for large $|x|$.  In this case,  for  $1<\alpha  \leq \alpha_F(\gamma,d,a),$  Pinsky \cite[p.153]{Pinsky1997}  proved \eqref{heat equation} does not posses global solution for any $u_0>0,$ and establish  classical global solutions  for $\alpha>  \alpha_{F}(\gamma, d).$ See \cite[p.153]{Pinsky1997}, \cite[Theorem 1]{Pinsky2009}.
\item \label{hr6}  Let  $d\geq 3,$ $\alpha= \alpha_F(d,0)$ or $1$,  and  $ V(x)=\frac{a}{1+|x|^b} \ (b>0)$   in \eqref{d3}.    In this case,  Zhang \cite{Zhang2001} found Fujita exponents 
under certain conditions on $a,b$. Later,  Ishige \cite[Theorems 1.1, 1.2]{Ishige2007}  considered $d\geq 2,$ and potential $V(x)= \frac{a}{|x|^2}$ with $a>0,$ and $b=1$ and determined the Fujita exponent $\alpha_F(d,0,a)$.  
See also recent work of Ishige and Kawakami in  \cite{Ishige2020}.
\end{enumerate}

\end{enumerate}
\end{remark}
\subsection{Dynamics of  HHE with inverse square potential}
We are now ready to state our well-posedness result in the following theorem.

\begin{theorem}[Well-posedness: subcritical and critical case]\label{lwp0}
Let $ q\in(1,\infty)$ and $\sigma_-,\sigma_+$ be as defined in \eqref{sigma}. Let 
\begin{equation}\label{E1}
\gamma\in\begin{cases}
(-2,\infty) & \text{if} \  a\leq 0\\
\R & \text{if} \  a> 0
\end{cases}
\end{equation}
and $\alpha$ satisfies
\begin{eqnarray}\label{E2}
\alpha\in \begin{cases} 
\ \quad \left(1
,1+\frac{\gamma+2}{\sigma_-}\right) & \text{if} \  a\leq 0\\
\left(1+\max\left(\frac{\gamma+2}{\sigma_-},
0\right),\infty\right) & \text{if} \   a>0
\end{cases}.
\end{eqnarray}
Let $s\ge \frac{\gamma}{\alpha-1}$, $s>\sigma_--\frac{d}{\alpha}$ and $\tau,\tau_c$ be as in \eqref{T} and satisfy
\begin{equation}\label{f02}
 \sigma_-<\tau<\sigma_++2\qquad
\text{and}\qquad \tau\le\tau_c.
\end{equation}
Then Cauchy problem \eqref{heat equation} is locally well-posed in $L_s^q(\R^d),$  and for the critical case we also have small data global existence.
	In the subcritical case,  if  we impose further restriction \begin{equation}\label{f01}
q>\alpha\qquad\text{and}\qquad \tau< \frac{\sigma_++2+\gamma}{\alpha},
\end{equation}then one has unconditional uniqueness\footnote{uniqueness is said to be unconditional if one can choose $X_T=C([0,T),Y)$ in the Definition \ref{well} for well-posedness so that  there is no need to find proper subset to achieve uniqueness}  in $C([0,T_m),L_s^q(\rd))$.
\end{theorem}

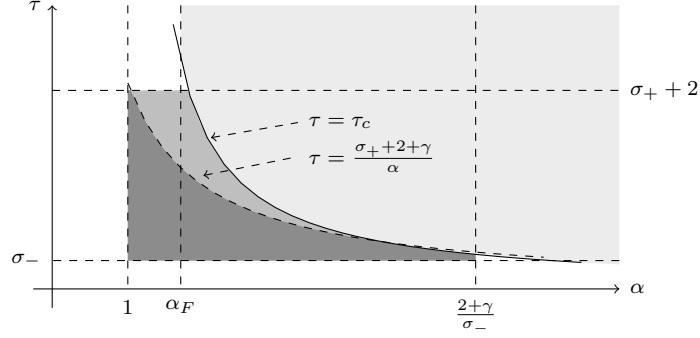
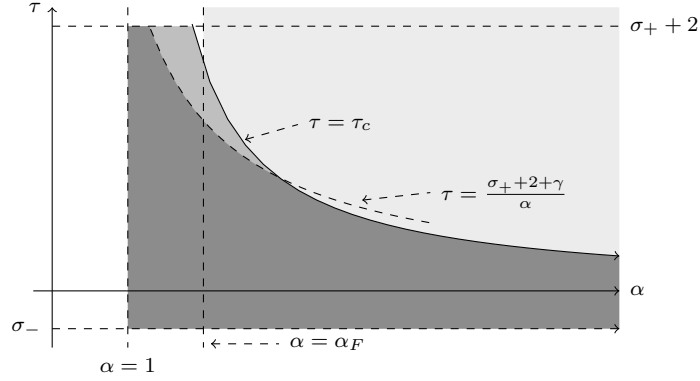
\begin{figure}
\subfigure[\tiny{The case  $d=3$, $a=-\frac{15}{64}$, $\gamma=\frac1{10}$ and $s\ge\frac{1}{\alpha-1}, s>\frac{3}{8}-\frac{3}{\alpha}$}] {
\begin{tikzpicture}[scale=1]

\fill [gray!50!white](1,3/8)--(1,21/8)--(1.8,21/8)--plot[domain=1.8:5.6]  ({\x},{(21/10)/(\x-1)})--(5.6,3/8);
\fill [gray!90!white](1,3/8)--(1,21/8)--(1+4/105,21/8)--plot[domain=1+4/105:109/25]  ({\x},{(109)/(40*\x)})--plot[domain=109/25:5.6]  ({\x},{(21/10)/(\x-1)})--(5.6,3/8);
\fill [gray!15!white](1.7,3)--(1.7,3.75)--(7.5,3.75)--plot[domain=7.5:1.7]  ({\x},{2.1/(\x-1)});
\draw[][->] (-.25,0)--(7.5,0) node[anchor=west] {\tiny{$\alpha$}};
\draw[][->] (0,-.25)--(0,3.75) node[anchor=east] {\tiny{$\tau$}};
\draw[dashed](1,3.75)--(1,0) node[anchor=north] {\tiny{$1$}};
\draw[dashed](5.6,3.75)--(5.6,0) node[anchor=north] {\tiny{$\frac{2+\gamma}{\sigma_-}$}};
\draw[][-] [domain=7:1.6] plot ({\x}, {(21/10)*(1/(\x-1))});
\draw[dashed][] [domain=1:6.5] plot ({\x}, {(109)/(40*\x)});
\draw[dashed][<-] (2.1,2)--(3.25,2.2) node[anchor= west] {\tiny{{$\tau=\tau_c$}}};
\draw[dashed][<-] (2,1.5)--(3.25,1.75) node[anchor= west] {\tiny{{$\tau=\frac{\sigma_++2+\gamma}{\alpha}$}}};
\draw[dashed](7.5,3/8)--(0,3/8) node[anchor=east] {\tiny{$\sigma_-$}};
\draw[dashed](0,21/8)--(7.5,21/8) node[anchor=west] {\tiny{$\sigma_++2$}};
\draw[dashed][-](1.7,3.75)--(1.7,0) node[anchor=north] {\tiny{$\alpha_F$}};
\end{tikzpicture}
}\\
\subfigure[\tiny{The case $d=3$, $a=\frac{3}{4}$, $\gamma=1$  and $s\ge\frac{1}{\alpha-1},s> -\frac{1}{2}-\frac{3}{\alpha}$}]  {
\begin{tikzpicture}[scale=1]
\fill [gray!50!white](1,-.5)--(1,3.5)--plot[domain=13/7:7]({\x},{3/(\x-1)})--(7,-.5);
\fill [gray!90!white](1,-.5)--(1,3.5)--plot[domain=9/7:3]({\x},{9/(2*\x)})--plot[domain=3:7.5]({\x},{3/(\x-1)})--(7.5,-.5);
\fill [gray!15!white](2,2.5)--(2,3.75)--(7.5,3.75)--plot[domain=7.5:2]  ({\x},{3/(\x-1)});
\draw[dashed](1,3.75)--(1,-.75) node[anchor=north] {\tiny{$\alpha=1$}};
\draw[][->] (-.25,0)--(7.5,0) node[anchor=west] {\tiny{$\alpha$}};
\draw[][->] (0,-.75)--(0,3.75) node[anchor=east] {\tiny{$\tau$}};
\draw[dashed][] [domain=1.3:5] plot ({\x}, {(9)/(2*\x)});
\draw[][<-] [domain=7.5:1.85] plot ({\x}, {(3*(1/(\x-1))});
\draw[dashed][<-] (2.55,2)--(3.25,2.2) node[anchor= west] {\tiny{{$\tau=\tau_c$}}};
\draw[dashed][<-] (4.1,1.2)--(5,1.3) node[anchor= west] {\tiny{{$\tau=\frac{\sigma_++2+\gamma}{\alpha}$}}};
\draw[dashed][<-](7.5,-.5)--(0,-.5) node[anchor=east] {\tiny{$\sigma_-$}};
\draw[dashed](0,3.5)--(7.5,3.5) node[anchor=west] {\tiny{$\sigma_++2$}};
\draw[dashed][-](2,3.75)--(2,-.75);
\draw[dashed][<-] (2.1,-.7)--(3,-.7) node[anchor= west] {\tiny{{$\alpha=\alpha_F$}}};
\end{tikzpicture}
}
\caption{\tiny{Local well-posedness in  $L_s^q(\rd)$ occurs in the  deep \& medium dark region by  Theorem \ref{lwp0} (only the boundary $\tau=\tau_c$ is included). Unconditional uniqueness in $L_s^q(\rd)$ is guaranteed by  Theorem \ref{lwp0} (furthermore part)  in the open deep dark region if $q>\alpha$. 
No LWP 
 in the unbounded lightest region
  by Theorem \ref{t:nonex}.}}  \label{f2}
\end{figure}

Theorem \ref{lwp0} is new for $a\neq 0$ and $\gamma>0.$ Up to now,  we could not know the well-posedness of \eqref{heat equation} with $\gamma>0$ in the mere $L^q-$spaces but in weighted $L^q_s-$spaces.  See Remark \ref{HR}\eqref{HRD}.  We prove  Theorem \ref{lwp0} via fixed point argument.  To this end,  the main new ingredient required is  our fixed-time estimate established in Theorem \ref{wle}.  

\begin{remark}\label{wgr} We have  several comments on Theorem \ref{lwp0}.
\begin{enumerate}

\item[--]  Theorem \ref{lwp0} recover results mentioned in  Remark \ref{hr}\eqref{hr1} and  is the main part  of a  detailed well-posedness  Theorem \ref{lwp2}. 
  
\item[--] For  $a=0$ and  $\tau=\tau_c$,  we  have from \eqref{f02} that   $\tau_c<d\Longleftrightarrow\alpha> \alpha_F$.  In this case, Theorem \ref{lwp0} 
along with below Theorem \ref{lwp2},    recover \cite[Theorem 1.4]{chikami2022optimal} and remove the assumption $q\ge\alpha$ and allows $s=\frac{\gamma}{\alpha-1}$.   See Remark \ref{hr}\eqref{hr4}. 

\item[--] {In \cite[Appendices B, C]{tayachi2023new}, equation \eqref{heat equation} is considered  for $a=0$,  to study local-wellposedness and lifespan of solutions in $L_s^q(\rd)$, $s\ge0$. However, these results do not cover the case when $s<0$ or $\tau=\tau_c$ or $\gamma\le-\min(2,d)$. Theorem \ref{lwp0} extends these local well-posedness results, particularly Theorems B1, C1 therein.  
}
\item[--] For $a=0,$ Theorem \ref{lwp0} eliminate technical hypothesis  (1.13) and 
$\alpha>\alpha_F$ from \cite[Theorem 1.13]{chikami2022optimal} in the subcritical case.

\item[--] Assume $s=0,\gamma<0$. Then for $a=0$ Theorem \ref{lwp0} recovers \cite[Theorem 1.1]{slimene2017well} and for $a\neq0$  Theorem \ref{lwp0} recovers \cite[Theorem 1.1]{Divyang_Saikatul_arxiv}. 

\item[--] For $V(x)= \frac{a_*}{|x|^2}$ and $d\geq 3$ in  \eqref{d3},  Ioku and Ogawa  \cite[Theorem 1.4]{Ioku2019} proved small data global existence for $1+ \frac{4}{d+2} < \alpha < 1+ \frac{4}{d-2}$. Theorem \ref{lwp0} relaxes this assumption and prove the result for any $\alpha>\alpha_F$ (note that $\alpha_F<1+ \frac{4}{d+2}$ for $d\ge2$).   See Remark \ref{rfbl}.

\item[--] In the subcritical case with  assumption \eqref{f01},
 Theorem \ref{lwp0} shows unconditional uniqueness of solution in $C([0,T_m),L_s^q(\rd))$ which complements \cite[Theorem 1.2]{chiiketanitay_arXiv}.  On the other hand, \cite[Theorem 1.13]{chikami2022optimal} established  uniqueness  for \eqref{d2} in a  proper subset of $C([0,T_m),L_s^q(\rd))$.
 See Remark \ref{rR}.   

\item[--] For detail comments on hypotheses of Theorem   \ref{lwp0},   see Remarks \ref{hy3},  \ref{hpcp}, \ref{fr}.
\end{enumerate}
\end{remark}

%

We now strengthen and complement  Theorem \ref{lwp0} by establishing following result.

\begin{theorem}[Finite time blow-up  for large data in the  subcritical case]
\label{blow-up large data}  Assume that $\tau\le \tau_c.$
Let $d,\gamma, \alpha,q,s $ be as in Theorem \ref{lwp0} (so local wellposedness for \eqref{heat equation} holds).  Let $F_\alpha$ satisfies  $F_\alpha(z)=z^\alpha$ for $z\geq0$\footnote{for example $F_\alpha(z)=\mu|z|^{\alpha-1}z, \mu|z|^\alpha$ or $\mu z^\alpha$}. 
Further assume
\begin{equation}\label{scon0} d+\gamma<
\begin{cases}\quad\alpha d &\text{ if }a=0\\
\alpha (d-2) &\text{ if }a\neq0
\end{cases}.
\end{equation}
Then there exists initial data $u_0\in L_s^q(\rd)$ such that $T_m(u_0)<\infty$. Moreover if $\tau<\tau_c$, one has a unique blow-up solution to \eqref{heat equation} with initial data $u_0$ in the following sense: 
there exist a unique solution $u$ of \eqref{heat equation} defined on $[0, T_m)$ such that 
\[ T_m< \infty \quad \text{and} \quad  \lim_{t\uparrow T_m}\|u(t)\|_{L_s^q}=\infty.\]
\end{theorem}

\begin{remark}\label{rfbl} We have  several comments for Theorem \ref{blow-up large data}.
\begin{itemize}
\item[--] For the critical case $\tau=\tau_c$,  similar blowup happens in a Kato  norm: If $T_m<\infty$,  one would have  $\|u\|_{\mathcal{K}_{k,s}^{p,q}(T_m)}=\infty$ for certain choice of $(k,p)$. See Section \ref{S4} for definition of Kato norm.
\item[--] Take $\gamma=s=0$ in \eqref{T},  and so   $\tau <\tau_c\Leftrightarrow q> \frac{d(\alpha-1)}{2}$.   Weissler \cite{weissler81} established blow-up solution for \eqref{d1}  in $L^q(\mathbb R^d).$  Theorem \ref{blow-up large data}  is compatible with this classical result.
\item[--] For $V(x)= \frac{a_*}{|x|^2},$  $u_0\in L^{\frac{d(\alpha-1)}{2}}(\mathbb R^d)$ with $\alpha\leq 1+ \frac{4}{d+2}$,   Ioku and Ogawa   \cite{Ioku2019}  pointed out that  \eqref{d3} have  blow-up solution  in finite time in $L^{\infty}-$norm.   However,  we are not aware of  any previous  results on finite time blow-up  solution in  $L^q_s-$norm   for $a, s, \gamma \neq 0$ and $q\neq \infty.$ 
 Thus Theorem \ref{blow-up large data}  is new.   
\item[--] 
Assume $d\ge3$, 
\begin{equation*}
\begin{cases} 1+\frac{\gamma+2}{d-2}<\alpha<1+\frac{\gamma+2}{\sigma_-}\ \quad &\text{for} \ a\le0\\
 1+\frac{\gamma}{d}<\alpha<\infty\ \quad &\text{for} \ a=0\\
1+\max(\frac{\gamma+2}{\sigma_-},\frac{\gamma+2}{d-2})<\alpha<\infty\quad &\text{for}\ a>0
\end{cases}
\end{equation*}
and  the hypothesis on $\gamma,q,s$ from Theorem \ref{lwp0}. Let $F(z)= |z|^\alpha$ or $|z|^{\alpha-1}z$ or $z^\alpha$. 
 Then  Theorem \ref{blow-up large data} reveals that,  
 there exists data in $L_s^q(\rd)$ such that the  local solution  established in Theorem \ref{lwp0} cannot be extend to global in time.  In the critical case, it also says that small data assumption in Theorem \ref{lwp0} is essentially optimal to establish global existence.

\end{itemize}
\end{remark}

\begin{definition}[weak solution] \label{weakdef}
 Let $u_0\in L^1_{loc}(\R^d)$, then we say a function $u$ is a weak solution to \eqref{heat equation} 
if $u\in L^{\alpha}( (0,T),  (L^{\alpha}_{\frac{\gamma}{\alpha}})_{loc}(\R^d))$ 
and satisfies the equation \eqref{heat equation} in the distributional sense, i.e.
\begin{align}\label{weak}
\notag\int_{\R^d} &u(T',x) \eta (T',x) \, dx-\int_{\R^d} u_0(x) \eta (0,x) \, dx\\
	&= \int_{[0,T']\times\R^d} u(t ,x)(\partial_t\eta+\Delta \eta-a|x|^{-2}\eta) (t ,x) 
	+ |x|^{\gamma} F_\alpha(u(t, x)) \,\eta(t,x)  \, dx\,dt
\end{align}
for all $T'\in [0,T]$ and for all $\eta \in C^{1,2}([0,T]\times \R^d)$ such that 
$\operatorname{supp} \eta(t, \cdot)$ is compact. The time $T$ is said to be the maximal existence time, which is denoted 
by $T_m^w$, if the weak solution cannot be extended beyond $[0,T).$
\end{definition}

\begin{remark}
Proceeding as \cite[Proposition 3.1]{ikeda2013small} it follows that $L_s^q$-integral solutions are weak solution. In that case $T_m\le T_m^w$.
\end{remark}

We shall now turn our attention to supercritical case.  In this case,   we show that  there exists positive initial data in $L^q_s(\rd)$ that do not generate a (weak) local solution to  \eqref{heat equation}. 
 Specifically,  we have the following theorem.
\begin{theorem}[Nonexistence of local positive weak solution in supercritical case]
\label{t:nonex}
Let $d\in \mathbb N$, $a,\gamma \in \mathbb R$, $\alpha$ satisfy \eqref{scon0} and
\begin{equation}\label{scon}
 \alpha>\alpha_F(d,\gamma,0)=1+\frac{(2+\gamma)^+}{d}.
\end{equation}
Assume that  $F_\alpha$ satisfies $F_\alpha(z)=z^\alpha$ for $z\geq0$,
 $q\in [1,\infty]$, $s\in\R.$ Let $\tau,\tau_c$ be as in \eqref{T}  and satisfy  $\tau< \tau_c.$
Then there exists an initial data $u_0 \in L^q_s (\R^d)$ such that  \eqref{heat equation} with $u(0)=u_0$ has no positive local weak solution. 
\end{theorem}
\begin{remark}\label{npsr} \
\begin{itemize}
\item[--]  
For $a=0=\gamma,$ Theorem \ref{t:nonex} recovers results mentioned in Remark \ref{hr}\eqref{hr2}.
\item[--] For $a=0,\gamma>-2$,   condition  $\alpha> \alpha_F(d,\gamma,0)$ \eqref{scon} implies  $d+\gamma< \alpha d$ in  \eqref{scon0}.  Thus,  in this case,   Theorem \ref{t:nonex} recovers  \cite[Theorem 1.16]{chikami2022optimal}. 
\item[--] Theorem \ref{t:nonex} implies failure of LWP in super-critical case. Theorem  \ref{t:nonex}  tells  if $\alpha$ satisfies \eqref{scon} then the sub criticality or criticality condition is necessary in Theorem \ref{lwp0}.   
\end{itemize}
\end{remark}
The paper is organized as follows.  In Section \ref{s2}, we gather some general tools which will be used later.  In Section \ref{wlep},   we prove Theorem \ref{wle}. In Section \ref{S4} we establish wellposedness  results.  In  Section \ref{S5},  we prove Theorems  \ref{blow-up large data} and \ref{t:nonex}.
\section{Preliminaries}\label{s2}
\noindent
\textbf{Notations}: The symbol $\alpha\wedge \beta$ means $\min(\alpha,\beta)$ whereas $\alpha\vee \beta$ mean $\max(\alpha,\beta)$. By $a^+$ we denote $a\vee0$.
The notation $A \lesssim B $ means $A \leq cB$ for some universal constant $c > 0.$ By $A\gtrsim B$ we mean $B\lesssim A$. By $A\sim B $ we mean $A\lesssim B$ and $A\gtrsim B$. 

  We shortly denote unweighed Lebesgue space norm by  $\|f\|_{L^p}= \|f\|_{p}.$ The Schwartz space is denoted by  $\mathcal{S}(\mathbb R^{d})$,  and the space of tempered distributions is  denoted by $\mathcal{S'}(\mathbb R^{d}).$ 
For $s\in \R$ and $q\in [1,\infty]$, we introduce the weighted local Lebesgue space $L^q_{s,loc}(\R^d)$ given by
\[L^q_{s,loc}(\R^d):=\left\{ f \in L^0(\R^d):
	f|_K\in L^q_s(\R^d),\ \forall K\subset \R^d,\ K\ \text{compact}
	\right\}
\] where $L^0(\rd)$ is the set of measurable functions on $\rd$.

\subsection{Lorentz space}
The Lorentz space is the space of all complex-valued measurable functions $f$ such that $\| f\|_{L^{p,q}(\rd)}<\infty$ where $\| f\|_{L^{p,q}(\rd)}$ is defined by 
\begin{equation*}\label{ls}
\| f\|_{L^{p,q}(\rd)}:=p^{\frac{1}{q}}\left\| t\mu\{|f|>t\}^{\frac{1}{p}}\right\|_{L^q\left((0,\infty),\frac{dt}{t}\right)}
\end{equation*} with $0<p<\infty$, $0<q\leq\infty$ and $\mu$ denotes the Lebesgue measure on $\rd$.  For shorthand,  we write
\begin{equation*}
\| f\|_{p,q}:=\| f\|_{L^{p,q}(\rd)}=\begin{cases}
p^{1/q}\left(\int_0^\infty t^{q-1}\mu\{|f|>t\}^{\frac{q}{p}}dt\right)^{1/q}\quad \text{for }q<\infty\\
\sup_{t>0}t\mu\{|f|>t\}^{\frac{1}{p}}\quad\quad\quad \text{for }q=\infty.
\end{cases}
\end{equation*}
Let us gather some useful results on Lorentz spaces relevant to subsequent  our proofs.
\begin{lemma}[Lemmata 2.2, 2.5 in \cite{o1963convolution}]Let $1\le p\le\infty$, $1\le q_1,q_2\le\infty$. Then
\begin{enumerate}
\item $\| f\|_{p,p}\sim\| f\|_{p}$, the usual Lebesgue $p$-norm.
\item $\| f\|_{p,q_2}\lesssim\| f\|_{p,q_1}$ if $q_1\ge q_2$.
\item $|\cdot|^{-b}\in L^{\frac{d}{b},\infty}(\rd)$ for $b>0$.
\end{enumerate}
\end{lemma}
\begin{lemma}[Theorems 2.6, 3.4 in \cite{o1963convolution}]\label{L}
We have the following inequalities in Lorentz spaces:
\begin{enumerate}
\item (H\"older's inequality) Let $\frac{1}{r}=\frac{1}{r_0}+\frac{1}{r_1}\in[0,1)$ and $s\geq1$ is such that $\frac{1}{s}\leq\frac{1}{s_0}+\frac{1}{s_1}$. Then $\| fg\|_{L^{r,s}}\leq r'\| f\|_{L^{r_0,s_0}}\| g\|_{L^{r_1,s_1}}$.
\item (Young's inequality) Let $\frac{1}{r}=\frac{1}{r_0}+\frac{1}{r_1}-1\in(0,1]$ and $s\geq1$ is such that $\frac{1}{s}\leq\frac{1}{s_0}+\frac{1}{s_1}$. Then $\| f\ast g\|_{L^{r,s}}\leq 3r\| f\|_{L^{r_0,s_0}}\| g\|_{L^{r_1,s_1}}$.
\end{enumerate}
\end{lemma}
 \subsection{Heat kernel estimate}\label{HK}
Let $g_a$ be the symmetric (in $x,y$ variable) heat kernel associated with the operator $\la$, i.e. 
\[ e^{-t \mathcal{L}_a} f(x)= \int_{\mathbb R^d } g_a(t,x, y) f(y) dy \quad  (t>0)\]see \cite[Proposition 3.6.]{metafune2021maximal}.
Then we have the following bounds for $g_a$: 

\begin{itheorem}[see Theorem 6.2 in \cite{metafune2018sharp}]\label{heat kernel}  Let $\sigma_-,\sigma_+$ be as defined in \eqref{sigma}. 
Let $d\ge2$, $a\geq a_*$.
Then there exist $c_1,c_2>0$ such that for any $t>0$ and $x,y\in\rd\backslash\{0\}$, the following estimate holds: 
\begin{equation*}
\big(1\vee\frac{\sqrt{t}}{|x|}\big)^{\sigma_-}\big(1\vee\frac{\sqrt{t}}{|y|}\big)^{\sigma_-} t^{-\frac{d}{2}}e^{-\frac{|x-y|^2}{c_1t}}\lesssim\ g_a(t,x,y)\ \lesssim \big(1\vee\frac{\sqrt{t}}{|x|}\big)^{\sigma_-}\big(1\vee\frac{\sqrt{t}}{|y|}\big)^{\sigma_-} t^{-\frac{d}{2}}e^{-\frac{|x-y|^2}{c_2t}}.
\end{equation*}
\end{itheorem}

\section{Dissipative estimates in weighted  Lebesgue spaces}\label{wlep}

  In order to prove Theorem \ref{wle},  we  first show it is enough to prove for $t=1$ (Lemma 3.1), then  using a duality (Lemma 3.2) we show it  is enough to prove for $s_1\ge0$. Then we crucially use  a known heat kernel estimate  (Theorem  \ref{heat kernel}) to  achieve the desired result.

\begin{lemma}\label{lem:scaling}
Let $1\le q_1,q_2 \le \infty$, and $s_1,s_2\in \mathbb R.$  Then 
$e^{-\mathcal{L}_a}$ is bounded from $L^{q_1}_{s_1}(\mathbb R^d)$ into $L^{q_2}_{s_2}(\mathbb R^d)$ if and only if 
$e^{-t\mathcal{L}_a}$ is bounded from $L^{q_1}_{s_1}(\mathbb R^d)$ into $L^{q_2}_{s_2}(\mathbb R^d)$ with 
\begin{equation}\label{linear-weight1'}
\|e^{-t\mathcal{L}_a}\|_{L^{q_1}_{s_1} \to L^{q_2}_{s_2}} =  t ^{-\frac{d}{2} (\frac{1}{q_1} - \frac{1}{q_2}) - \frac{s_1 - s_2}{2}} \|e^{-\mathcal{L}_a}\|_{L^{q_1}_{s_1} \to L^{q_2}_{s_2}}
\end{equation}
for any $t>0$. 
\end{lemma}
\begin{proof} 
It is enough to show \eqref{linear-weight1'} if $e^{-\mathcal{L}_a}$ is bounded from $L^{q_1}_{s_1}(\mathbb R^d)$ into $L^{q_2}_{s_2}(\mathbb R^d)$, since the converse is trivial. 
The proof is based on the scaling argument. 
Let  $f \in L^{q_1}_{s_1}(\mathbb R^d)$. 
Since
\[
(e^{-t\mathcal{L}_a} f)(x) =
\left(e^{-\mathcal{L}_a} (f(t^{\frac12} \cdot))\right)(t^{-\frac12} x),
\]
\[
(e^{-\mathcal{L}_a} f)(x) =
\left(e^{-t\mathcal{L}_a} (f(t^{-\frac12} \cdot))\right)(t^{\frac12} x),
\]
for $t>0$ and $x\in \mathbb R^d$, we have
\[
\|e^{-t\mathcal{L}_a} f\|_{L^{q_2}_{s_2}}
\le  t ^{-\frac{d}{2} (\frac{1}{q_1} - \frac{1}{q_2}) - \frac{s_1 - s_2}{2}} \|e^{-\mathcal{L}_a}\|_{L^{q_1}_{s_1} \to L^{q_2}_{s_2}}
\|f\|_{L^{q_1}_{s_1}},
\]
\[
\|e^{-\mathcal{L}_a} f\|_{L^{q_2}_{s_2}}
\le  t ^{\frac{d}{2} (\frac{1}{q_1} - \frac{1}{q_2}) + \frac{s_1 - s_2}{2}}
\|e^{-t\mathcal{L}_a}\|_{L^{q_1}_{s_1} \to L^{q_2}_{s_2}}
\|f\|_{L^{q_1}_{s_1}}.
\]
Hence, \eqref{linear-weight1'} is proved. 
\end{proof}

\begin{lemma}\label{duality}
Let $q_1,q_2\in(1,\infty)$ and $s_1,s_2\in\R$ and $A=\{x\in\rd:|x|\ge1\}$. 
Let $k(x,y)=k(y,x)$ for $x,y\in A$ and for $x\in A$ set
$Tf(x)=\int_A k(x,y)f(y)dy$.
Then \begin{equation*}
   \|Tf\|_{L^{q_2}_{s_2}(A)}\le C\|f\|_{L^{q_1}_{s_1}(A)} \text{ for all }f
\end{equation*}
if and only if
\begin{equation*}
   \|Tf\|_{L^{q_1'}_{-s_1}(A)}\le C\|f\|_{L^{q_2'}_{-s_2}(A)} \text{ for all }f.
\end{equation*}
\end{lemma}
\begin{proof}
Note that 
\begin{eqnarray*}
 \|Tf\|_{L^{q_1'}_{-s_1}(A)}
  &=&\sup_{\|g\|_{L_{s_1}^{q_1}}\le1}|\int_{A}\int_{A}k(x,y)f(y)dyg(x)dx|\\
    & = & \sup_{\|g\|_{L_{s_1}^{q_1}}\le1}|\int_{A}\int_{A}k(x,y)g(x)dxf(y)dy |\\
    &=&\sup_{\|g\|_{L_{s_1}^{q_1}}\le1}|\int_{A}(Tg)(y)f(y)dy | \le \sup_{\|g\|_{L_{s_1}^{q_1}}\le1}\|Tg\|_{L_{s_2}^{q_2}}\|f\|_{L_{-s_2}^{q_2'}}\\
    &\le&c\sup_{\|g\|_{L_{s_1}^{q_1}}\le1}\|g\|_{L_{s_1}^{q_1}}\|f\|_{L_{-s_2}^{q_2'}}=c\|f\|_{L_{-s_2}^{q_2'} (A)}
\end{eqnarray*}
This completes the proof.
\end{proof}
\begin{proof}[\textbf{Proof of Theorem \ref{wle} (Sufficiency part)}] Assume that \eqref{Condiexp} and \eqref{reg} hold.
In view of Lemma \ref{lem:scaling} it is enough to prove the case $t=1$ i.e.
\begin{equation}\label{ep0}
 \|e^{-\mathcal{L}_a}f\|_{L^{q_2}_{s_2}}\lesssim\|f\|_{L^{q_1}_{s_1}}.
\end{equation} 
For $x\in\rd$ and $f\in L_{s_1}^{q_1}(\rd)$ applying Theorem \ref{heat kernel} we achieve
\begin{equation}\label{ep1}
|e^{- \mathcal{L}_a}f(x)|
\lesssim\ \left(1\vee\frac{1}{|x|}\right)^{\sigma_-}
\int_{\R^d}\left(1\vee\frac{1}{|y|}\right)^{\sigma_-}G(x-y)|f(y)|dy.
\end{equation} where 
    $G(x):=e^{-\frac{|x|^2}{c_2}}$%
with $c_2$ as in Theorem \ref{heat kernel}. 
Set
\[
1_{\ge 1}(x):=\begin{cases}
0\text{ for }|x|< 1\\
1 \text{ for }|x|\ge 1
\end{cases}\qquad\text{and}\qquad1_{<1}:=1-1_{\ge1}.
\] 
Then  using $e^{-\mathcal{L}_a}f=1_{\ge1}e^{-\mathcal{L}_a}f+1_{<1}e^{-\mathcal{L}_a}f$ and \eqref{ep1} we have
\begin{eqnarray*}
 \|e^{-\mathcal{L}_a}f\|_{L^{q_2}_{s_2}} & \le&  \|e^{-\mathcal{L}_a}f(x)\|_{L^{q_2}_{s_2}{(|x|\ge 1)}}+\|e^{-\mathcal{L}_a}f(x)\|_{L^{q_2}_{s_2}{(|x|<1)}}\\
 & \lesssim&  \|\int_{\R^d}\big(1\vee\frac{1}{|y|}\big)^{\sigma_-}G(x-y)|f(y)|dy\|_{L^{q_2}_{s_2}{(|x|\ge 1)}}\\
 &&+\ \| |x|^{-\sigma_-}
\int_{\R^d}\big(1\vee\frac{1}{|y|}\big)^{\sigma_-}G(x-y)|f(y)|dy\|_{L^{q_2}_{s_2}{(|x|<1)}}.
\end{eqnarray*}
Splitting the integrations  in $y$ variable we obtain 
\begin{eqnarray*}
\|e^{-\mathcal{L}_a}f\|_{L^{q_2}_{s_2}} & \lesssim&  \|\int_{|y|\ge1}G(x-y)|f(y)|dy\|_{L^{q_2}_{s_2}{(|x|\ge 1)}}\\
&&+\ \|\int_{|y|<1}|y|^{-\sigma_-}G(x-y)|f(y)|dy\|_{L^{q_2}_{s_2}{(|x|\ge 1)}}\\
 &&+\ \| |x|^{-\sigma_-}
\int_{|y|\ge2}G(x-y)|f(y)|dy\|_{L^{q_2}_{s_2}{(|x|<1)}}\\
 &&+\ \| |x|^{-\sigma_-}
\int_{|y|<2}|y|^{-\sigma_-}G(x-y)|f(y)|dy\|_{L^{q_2}_{s_2}{(|x|<1)}}=:I+\textit{II}+\textit{III}+\textit{IV}.
\end{eqnarray*}

Now we show that each of these terms  is dominated by $ \|f\|_{L^{q_1}_{s_1}}$ which would prove \eqref{ep0} to conclude the proof.

{\it Estimate for $\textit{IV}$:} Using boundedness of $G$ and changing the order of integration and H\"older's inequality we obtain
\begin{eqnarray*}
\textit{IV}&\lesssim&\| |x|^{-\sigma_-}
\int_{|y|<2}|y|^{-\sigma_-}|f(y)|dy\|_{L^{q_2}_{s_2}(|x|<1)}\\
&=&\| |x|^{-\sigma_-}\|_{L^{q_2}_{s_2}(|x|<1)}
\int_{|y|<2}|y|^{-\sigma_-}|f(y)|dy\\
&=&\| |x|^{-\sigma_-}\|_{L^{q_2}_{s_2}(|x|<1)}
\int_{|y|<2}|y|^{-\sigma_--s_1}|y|^{s_1}|f(y)|dy\\
&\le&\| |x|^{-\sigma_-}\|_{L^{q_2}_{s_2}(|x|<1)}
\||y|^{-\sigma_--s_1}\|_{L^{q_1'}(|y|<2)}\||y|^{s_1}f(y)\|_{q_1}\ \lesssim\ \|f\|_{L^{q_1}_{s_1}},
\end{eqnarray*}
where in the last step we used the hypothesis \begin{eqnarray*}
(s_2-\sigma_-)q_2+d>0&\Longleftrightarrow& \sigma_-<s_2+\frac{d}{q_2},\\
(-\sigma_--s_1)q_1'+d>0&\Longleftrightarrow& s_1+\frac{d}{q_1}<d-\sigma_-=\sigma_++2.
\end{eqnarray*}

{\it Estimate for $\textit{III}$:} Note that for $|x|<1$,  $|y|\ge2$ we have  $|x-y|\ge |y|-|x|\ge \frac{1}{2}|y|$, and as $G$ is radially decreasing we have $G(x-y)\le G(\frac{y}{2}).$ Therefore,  we have 
\begin{eqnarray*}
\textit{III}&\le&\| |x|^{-\sigma_-}
\int_{|y|\ge2}G(\tfrac{y}{2})|f(y)|dy\|_{L^{q_2}_{s_2}{(|x|<1)}}\\
 &=& \| |x|^{-\sigma_-}\|_{L^{q_2}_{s_2}{(|x|<1)}}\int_{|y|\ge2}G(\tfrac{y}{2})|y|^{-s_1}|y|^{s_1}|f(y)|dy\\
 &\le& \| |x|^{-\sigma_-}\|_{L^{q_2}_{s_2}{(|x|<1)}}\|G(\tfrac{y}{2})|y|^{-s_1}\|_{L^{q_1'}(|y|\ge2)}\||y|^{s_1}f(y)\|_{q_1}\ \lesssim\ \|f\|_{L^{q_1}_{s_1}},
\end{eqnarray*}
where in the last step we  used the hypothesis $ \sigma_-<s_2+\frac{d}{q_2}$ as in the estimate for {\it IV} and the fact that $G$ is Schwartz class function. 

{\it Estimate for $\textit{II}$:}  We claim that $\| |x|^{s_2}G(x-y)\|_{L^{q_2}{(|x|\ge1)}}\lesssim 1$ uniformly for all $|y|<1$. In fact when $s_2\le0$ we have $\| |x|^{s_2}G(x-y)\|_{L^{q_2}{(|x|\ge1)}}\le \| G(x-y)\|_{L^{q_2}{(|x|\ge1)}}\le  \| G\|_{q_2}$ for all $y$. On the other hand when $s_2>0$, using $|x|^{s_2}\lesssim |x-y|^{s_2}+|y|^{s_2}$, for $|y|<1$ we have
\begin{eqnarray*}
\| |x|^{s_2}G(x-y)\|_{L^{q_2}{(|x|\ge1)}}&\lesssim &\| |x-y|^{s_2}G(x-y)\|_{L^{q_2}{(|x|\ge1)}}+\| |y|^{s_2}G(x-y)\|_{L^{q_2}{(|x|\ge1)}}\\
&\le &\| |\cdot|^{s_2}G\|_{q_2}+\|G(x-y)\|_{L^{q_2}{(|x|\ge1)}}\\
&\le &\| |\cdot|^{s_2}G\|_{q_2}+\|G\|_{q_2}.
\end{eqnarray*} 
This proves the claim.
Then 
\begin{eqnarray*}
\textit{II}&=&\||x|^{s_2}\int_{|y|<1}|y|^{-\sigma_-}G(x-y)|f(y)|dy\|_{L^{q_2}{(|x|\ge 1)}}\\
&\le&\|\int_{|y|<1}\||x|^{s_2}G(x-y)\|_{L^{q_2}{(|x|\ge 1)}}|y|^{-\sigma_-} |f(y)|dy\\
&\lesssim&\|\int_{|y|<1}|y|^{-\sigma_--s_1}|y|^{s_1} |f(y)|dy\\
&\le&\||y|^{-\sigma_--s_1}\|_{L^{q_1'}(|y|<1)}\||y|^{s_1}f(y)\|_{q_1}\ \lesssim\ \|f\|_{L^{q_1}_{s_1}}
\end{eqnarray*}
using the claim above and the hypothesis $s_1+\frac{d}{q_1}<\sigma_++2 $ as in the estimate for {\it IV}. 

{\it Estimate for $I$:} Let us treat $I$ case by case.  

\fbox{Case $s_1=s_2=0$} By hypothesis $\frac{1}{p}:=1+\frac{1}{q_2}-\frac{1}{q_1}\in[0,1]$. Then  using Young's inequality
\begin{eqnarray*}
I= \|\int_{|y|\ge1}G(x-y)|f(y)|dy\|_{L^{q_2}{(|x|\ge 1)}}\le\|G\ast (1_{\ge1}|f|)\|_{q_2}\le\|G\|_{p}\|f\|_{q_1}\lesssim\ \|f\|_{L^{q_1}_{0}}.
\end{eqnarray*}

\fbox{Case $0=s_2<s_1$} If $\frac{1}{q_2}=\frac{1}{q_1}+\frac{s_1}{d}$, then using Young's and Holder's inequalities in Lorentz spaces i.e. Lemma \ref{L} we have
\begin{eqnarray*}
I\ \lesssim\ \|G\|_{1,q_2}\|1_{\ge1}f\|_{q_2,\infty}\ \lesssim\ \||\cdot|^{-s_1}\|_{\frac{d}{s_1},\infty}\||\cdot|^{s_1}f\|_{q_1,\infty}\ \lesssim\ \|f\|_{L^{q_1}_{s_1}}.
\end{eqnarray*}
 If $\frac{1}{q_2}< \frac{1}{q_1}+\frac{s_1}{d}$, then by using Lemma \ref{tl} \eqref{tl1}  we choose $\frac{1}{p_0},\frac{1}{p_1},\frac{1}{p_2}\in[0,1]$ so that 
\begin{equation}\label{ep2}
1+\frac{1}{q_2}=\frac{1}{p_0}+\frac{1}{p_1} \qquad\frac{1}{p_1}=\frac{1}{p_2}+\frac{1}{q_1},\qquad \frac{1}{p_2}<\frac{s_1}{d}.
\end{equation}and then  using Young's and Holder's inequalities we achieve
\begin{eqnarray*}
I&\le& \|G\|_{p_0}\|1_{\ge1}f\|_{p_1}\\
&\lesssim&\|1_{\ge1}|\cdot|^{-s_1}\|_{p_3}\||\cdot|^{s_1}f\|_{q_1}\ \lesssim\ \|f\|_{L^{q_1}_{s_1}}.
\end{eqnarray*}

\fbox{Case $0<s_2=s_1$}
 Using $|x|^{s_2}\lesssim |x-y|^{s_2}+|y|^{s_2}$ and Young's inequality
\begin{eqnarray*}
I&=& \||x|^{s_2}\int_{|y|\ge1}G(x-y)|f(y)|dy\|_{L^{q_2}{(|x|\ge 1)}}\\
&\lesssim& \|\int_{|y|\ge1}|x-y|^{s_2}G(x-y)|f(y)|dy\|_{q_2}+\|\int_{|y|\ge1}G(x-y)|y|^{s_2}|f(y)|dy\|_{q_2}\\
&=&\|(|\cdot|^{s_2}G)\ast (1_{\ge1}|f|)\|_{q_2}+\|G\ast (1_{\ge1}|\cdot|^{s_2}|f|)\|_{q_2}:=Ia+Ib
\end{eqnarray*}
Note that we have $\frac{d}{q_2}< \frac{d}{q_1}+{s_1}$ as $s_2>0$ and \eqref{Condiexp} is assumed. Then by choosing $\frac{1}{p_0},\frac{1}{p_1},\frac{1}{p_2}\in[0,1]$ satisfying \eqref{ep2} and using Young's and Holder's inequalities we achieve
\begin{eqnarray*} 
Ia&\le& \||\cdot|^{s_2}G\|_{p_0}\|1_{\ge1}f\|_{p_1}\\
&\lesssim&\|1_{\ge1}|\cdot|^{-s_1}\|_{p_2}\||\cdot|^{s_1}f\|_{q_1}\ \lesssim\ \|f\|_{L^{q_1}_{s_1}}.
\end{eqnarray*}
By hypothesis  $\frac{1}{p}:=1+\frac{1}{q_2}-\frac{1}{q_2}\in[0,1]$ then 
\begin{eqnarray*}
Ib&\le&\|G\|_{p} \|1_{\ge1}|\cdot|^{s_2}f\|_{q_1}\lesssim\ \|f\|_{L^{q_1}_{s_1}}.
\end{eqnarray*}


\fbox{Case $0<s_2< s_1$} Since $\frac{d}{q_2}< \frac{d}{q_1}+{s_1}$ we proceed as in above case and prove the estimate for $Ia$. 
Now with the  assumption $s_2<s_1$  using Lemma \ref{tl} \eqref{tl1'}  we choose $\frac{1}{p_3},\frac{1}{p_4},\frac{1}{p_5}\in[0,1]$ satisfying
\begin{equation}\label{ep3}
1+\frac{1}{q_2}=\frac{1}{p_3}+\frac{1}{p_4},\qquad\frac{1}{p_4}=\frac{1}{p_5}+\frac{1}{q_2},\qquad\frac{1}{p_5}<\frac{s_1-s_2}{d}
\end{equation} and obtain
\begin{eqnarray*}
Ib&\le&\|G\|_{p_3} \|1_{\ge1}|\cdot|^{s_2}f\|_{p_4}\\
&\lesssim&\|1_{\ge1}|\cdot|^{s_2-s_1}\|_{p_5}\||\cdot|^{s_1}f\|_{q_1}\ \lesssim\ \|f\|_{L^{q_1}_{s_1}}.
\end{eqnarray*}

\fbox{Case $s_2<0< s_1$} If $\frac{d}{q_2}+s_2< \frac{d}{q_1}+{s_1}$,
by Lemma \ref{tl} \eqref{tl2} we choose $\frac{1}{p_6},\cdots,\frac{1}{p_{10}}\in[0,1]$ satisfying 
\begin{equation}\label{ep5}
\frac{1}{q_2}=\frac{1}{p_6}+\frac{1}{p_7},\qquad1+\frac{1}{p_7}=\frac{1}{p_8}+\frac{1}{p_9},\qquad\frac{1}{p_9}=\frac{1}{p_{10}}+\frac{1}{q_1}, \qquad \frac{1}{p_6}<-\frac{s_2}{d},\qquad\frac{1}{p_{10}}<\frac{s_1}{d}
\end{equation}
so that
\begin{eqnarray*}
I&=&\||x|^{s_2}\int_{|y|\ge1}G(x-y)|f(y)|dy\|_{L^{q_2}{(|x|\ge 1)}}\\
&\le&\||x|^{s_2}\|_{L^{p_6}(|x|\ge1)}\|G\ast (1_{\ge1}|f|)\|_{L^{p_7}{(|x|\ge 1)}}\\
&\le&\||x|^{s_2}\|_{L^{p_6}(|x|\ge1)}\|G\|_{p_8} \|1_{\ge1}f\|_{p_9}\\
&\le&\||x|^{s_2}\|_{L^{p_6}(|x|\ge1)}\|G\|_{p_8} \||y|^{-s_1}\|_{L^{p_{10}}(|y|\ge1)}\||\cdot|^{s_1}f\|_{q_1}\ \lesssim \|f\|_{L^{q_1}_{s_1}}.
\end{eqnarray*}
If $\frac{d}{q_2}+s_2=\frac{d}{q_1}+{s_1}$, then we claim $0<\frac{-s_2}{d}<1$. We need to show $-s_2<d$ i.e. $s_2>-d$.  Infact if $s_2\le- d$, then $\frac{d}{q_1}+{s_1}=\frac{d}{q_2}+s_2\le \frac{d}{q_2}-d<0$, a contradiction as $q_1,s_1>0$.
Next we claim $0<\frac{1}{q_1}+\frac{s_1}{d}<1$. This is  because $0<\frac{1}{q_1}+\frac{s_1}{d}=\frac{1}{q_2}+\frac{s_2}{d}<\frac{1}{q_2}<1$ using $s_2<0$.\\
Above claim shows $\frac{1}{p_{12}}:=\frac{s_1}{d}+\frac{1}{q_1}\in(0,1)$, then we have $\frac{1}{q_2}=\frac{1}{-d/s_2}+\frac{1}{p_{12}}$. Therefore
%
\begin{eqnarray*}
I&=&\||x|^{s_2}\int_{|y|\ge1}G(x-y)|f(y)|dy\|_{L^{q_2}{(|x|\ge 1)}}\\
&\le&\||x|^{s_2}\|_{\frac{d}{-s_2},\infty}\|G\ast (1_{\ge1}|f|)\|_{p_{12},q_2}\\
&\le&\||x|^{s_2}\|_{\frac{d}{-s_2},\infty}\|G\|_{1,q_{2}} \|1_{\ge1}f\|_{p_{12},\infty}\\
&\le&\||x|^{s_2}\|_{\frac{d}{-s_2},\infty}\|G\|_{1,q_{2}} \||y|^{-s_1}\|_{\frac{d}{s_1},\infty}\||\cdot|^{s_1}f\|_{q_1,\infty}\ \lesssim \|f\|_{L^{q_1}_{s_1}}.
\end{eqnarray*}


\fbox{Case $s_2\le s_1\le0$} Follows from  duality Lemma \ref{duality} and the above cases.
This completes the proof.
\end{proof}
\begin{proof}[\textbf{Proof of Theorem \ref{wle} (Necessity part)}] Assume that \eqref{main est} hold.
Let $g=\exp({-\frac{|\cdot|^2}{c_1}})$ where $c_1$ as in Theorem \ref{heat kernel}.

{\it Necessity of $\sigma_-<s_2+\frac{d}{q_2}$, $s_1+\frac{d}{q_1}<\sigma_++2$:}
Let $f$ be supported in $B(0,1)$ and equal to $|\cdot|^\theta$ in $B(0,\frac{1}{2})$ with $\theta>\max(-s_1-\frac{d}{q_1},\sigma_--d,0)$.  Then $f\in L_{s_1}^{q_1}(\rd)$ and hence by hypothesis \eqref{main est}, we have $e^{-\la }f\in L_{s_2}^{q_2}(\rd)$. On the other hand 
 for $|x|\le1$  
\begin{eqnarray*}
[e^{-\la }f](x)&\ge&\int_{|y|\le1/2}g_a(1,x,y)f(y)dy\\
&\gtrsim& |x|^{-\sigma_-}\int_{|y|\le1/2}|y|^{-\sigma_-+\theta} g(x-y)dy\sim |x|^{-\sigma_-}.
\end{eqnarray*}where we have used Theorem \ref{heat kernel} in the second step and $\theta>\sigma_--d$ in the last step. 
Since $e^{-\la }f\in L_{s_2}^{q_2}(\rd)$, we must  have $\sigma_-<s_2+\frac{d}{q_2}$. 
Using symmetry of heat kernel (see Remark \ref{HR} \eqref{symm}), it follows that $s_1+\frac{d}{q_1}<\sigma_++2$.  This proof is inspired by the ideas used in  \cite[Section 4]{killip2018sobolev} where $q_1=q_2$, $s_1=s_2=0$ was treated.

{\it Necessity of $s_2+\frac{d}{q_2}\le s_1+\frac{d}{q_1}$:} Let $ 0 \neq f\in L^2\cap L_{s_1}^{q_1}.$ If $s_2+\frac{d}{q_2}> s_1+\frac{d}{q_1}$, then using \eqref{main est}, we have $e^{-t\mathcal{L}_a}f\to0$ in $L^{q_2}_{s_2}$ (and hence pointwise a.e.) as $t\to0$. Since $f\in L^2$, using semigroup property,  we have $e^{-t\mathcal{L}_a}f\to f$ in $L^2$ as $t\to0$.
 Thus $f=0$ which is a contraction.

{\it Necessity of $s_2\le s_1$:} We prove this by modifying the proof in case $a=0$ in \cite[Remark 10]{tayachi2023new}. Let $\varphi\in L_{s_1}^{q_1}$ be a smooth non-negative function with support in $B(0,1)$ and take $f_\tau=\varphi(\cdot-\tau x_0)$ with $|x_0|=1$. Then, 
for $\tau>2$ and $|x|\ge1,$ we have    \begin{eqnarray*}
[e^{-\la}f_\tau](x)&\ge&\int_{|y|\ge1}g_a(1,x,y)f_\tau(y)dy\\
&\gtrsim&\int_{|y|\ge1} g(x-y)f_\tau(y)dy\\
&=&\int g(x-y)f_\tau(y)dy=(g\ast f_\tau) (x)=(g\ast \varphi)(\cdot-\tau x_0),
\end{eqnarray*}where we  used Theorem \ref{heat kernel} in the second step, the fact $B(0,1)\cap\supp  ( f_\tau)=\emptyset$ in the third step.
Now $\||\cdot|^{s_2}(g\ast \varphi)(\cdot-\tau x_0)\|_{q_2}=\||\cdot+\tau x_0|^{s_2}(g\ast \varphi)\|_{q_2}=\tau^{s_2}\||\frac{\cdot}{\tau}+ x_0|^{s_2}(g\ast \varphi)\|_{q_2}$ and $\||\cdot|^{s_1}f\|_{q_1}=\tau^{s_1}\||\frac{\cdot}{\tau}+ x_0|^{s_1}\varphi\|_{q_1}$. Therefore for $\tau>2$ we have from \eqref{main est} that
\[
 \tau^{s_2-s_1}\left\|\left|\frac{\cdot}{\tau}+ x_0\right|^{s_2}(g\ast \varphi)\right\|_{q_2}\lesssim\left\|\left|\frac{\cdot}{\tau}+ x_0\right|^{s_1`}\varphi\right\|_{q_1}
\]but $\||\frac{\cdot}{\tau}+ x_0|^{s_2}(g\ast \varphi)\|_{q_2}\to \|g\ast \varphi\|_{q_2}$ and $\||\frac{\cdot}{\tau}+ x_0|^{s_1`}\varphi\|_{q_1}\to\|\varphi\|_{q_1}$ as $\tau\to\infty$. Therefore we must have $s_2\le s_1$.
\end{proof}
\begin{lemma}\label{tl}
There exists $p_0,\cdots,p_{10}\in[1,\infty]$ so that
\begin{enumerate}
\item\label{tl1} if $0< s_1$, $\frac{d}{q_2}< \frac{d}{q_1}+{s_1}$ hold, then \eqref{ep2} is satisfied,
\item\label{tl1'} if $s_2< s_1$ holds, then   \eqref{ep3}  is satisfied,
\item\label{tl2} if $s_2<0< s_1$, $\frac{d}{q_2}+s_2< \frac{d}{q_1}+{s_1}$ hold, then  \eqref{ep5}  is satisfied.
\end{enumerate} 
\end{lemma}
\begin{proof}
\eqref{tl1} Choose $$\frac{1}{p_0}\in\big(\max(\frac{1}{q_2},1+\frac{1}{q_2}-\frac{1}{q_1}-\frac{s_1}{d}),\min(1,1+\frac{1}{q_2}-\frac{1}{q_1})\big).$$
 The last interval in nonempty as $q_1,q_2\in(1,\infty)$, $s_1>0$ and $\frac{d}{q_2}< \frac{d}{q_1}+{s_1}$.  Now set $
\frac{1}{p_1}=1+\frac{1}{q_2}-\frac{1}{p_0}$, $\frac{1}{p_2}=1+\frac{1}{q_2}-\frac{1}{p_0}-\frac{1}{q_1}$ then \eqref{ep2} is satisfied.

\eqref{tl1'} Proof is similar to \eqref{tl1}, only $s_1$ is replaced by $s_1-s_2$.
Choose $$\frac{1}{p_3}\in\big(\max(\frac{1}{q_2},1+\frac{1}{q_2}-\frac{1}{q_1}-\frac{s_1-s_2}{d}),\min(1,1+\frac{1}{q_2}-\frac{1}{q_1})\big).$$
 The last interval in nonempty as $q_1,q_2\in[1,\infty]$, $s_1-s_2>0$ and $\frac{1}{q_2}< \frac{1}{q_1}+\frac{s_1-s_2}{d}$. Now set $
\frac{1}{p_4}=1+\frac{1}{q_2}-\frac{1}{p_3}$, $\frac{1}{p_5}=1+\frac{1}{q_2}-\frac{1}{p_3}-\frac{1}{q_1}$ then \eqref{ep3} is satisfied.

\eqref{tl2} Note that $1+\frac{1}{q_2}+\frac{s_2}{d}-\frac{1}{q_1}-\frac{s_1}{d}<1$, then choose 
$$\frac{1}{p_8}\in(\max(1+\frac{1}{q_2}+\frac{s_2}{d}-\frac{1}{q_1}-\frac{s_1}{d},1-\frac{1}{q_1}-\frac{s_1}{d},\frac{1}{q_2}+\frac{s_2}{d},0),\min(1+\frac{1}{q_2}-\frac{1}{q_1},1)).$$ Then choose $$\frac{1}{p_7}\in(\max(\frac{1}{q_2}+\frac{s_2}{d},\frac{1}{p_8}+\frac{1}{q_1}-1,0),\min(\frac{1}{p_8}+\frac{1}{q_1}+\frac{s_1}{d}-1,\frac{1}{q_2},\frac{1}{p_8})).$$ Set
\[
\frac{1}{p_6}=\frac{1}{q_2}-\frac{1}{p_7},\qquad\frac{1}{p_9}=1+\frac{1}{p_7}-\frac{1}{p_8},\qquad\frac{1}{p_{10}}=1+\frac{1}{p_7}-\frac{1}{p_8}-\frac{1}{q_1}
\]
 so that equalities in \eqref{ep5} are satisfied.  
%
\end{proof}

\section{Local and small data global  well-posedness}\label{S4}

In this section we prove the well-posedness in critical and subcritical case i.e. when $\tau\le\tau_c$ (recall that $\tau=\frac{d}{q}+s$ and $\tau_c=\frac{2+\gamma}{\alpha-1}$). 
In order to prove Theorem \ref{lwp0}, we introduce the Kato space depending on four parameters $(p,q,k,s)$.

\begin{definition}[Kato space]\label{kato}
\label{def:Kato}
 Let $k,s\in \R$ and $p,q\in [1,\infty]$,  and set \[\beta=\beta(d,k,s,p,q):=\frac{1}{2}(s+\frac{d}{q}-k-\frac{d}{p}).\]
Then the Kato space $\mathcal{K}^{p,q}_{k,s}(T)$ is defined by 
\begin{equation}\nonumber
   \mathcal{K}^{p,q}_{k,s}(T)
   :=\left\{u:[0,T)\to L_k^p(\rd):
   	\|u\|_{\mathcal{K}^{p, q}_{k,s}(T')}
   <\infty\ \text{for any } T' \in (0,T)\right\}
\end{equation}
endowed with the norm
\[
\|u\|_{\mathcal K^{p,q}_{k,s}(T)}
	:=\sup_{0\le t\le T}t^{\beta} \|u(t)\|_{L^{p}_k}.
\]
\end{definition}
\begin{remark}\label{RR}
In \cite{chikami2022optimal}, Kato space with three parameter was used.  This is basically $\mathcal{K}^{q,q}_{k,s}(T)$ when one puts $p=q$ in Definition \ref{kato}. This restriction didn't allow  authors in \cite{chikami2022optimal} to consider the case $1\le q<\alpha$. See also \cite[Remarks 16, 18]{tayachi2023new}.
\end{remark}

By Theorem \ref{wle},  we immediately get   the following result (in fact these results are equivalent):
\begin{lemma}\label{Kle} Let $k,s\in \R$ and $p,q\in (1,\infty)$. Then  \[
   \|e^{-t\mathcal{L}_a}f\|_{\mathcal K^{p,q}_{k,s}(\infty)}\le C\|f\|_{L^q_s},\qquad\forall f\in L^q_s(\R^d)
\]if and only if \begin{equation}
\label{condition:p,k}
    k\le s\ \ \ \text{and}\ \ \sigma_-< \frac{d}{p}+k\le \frac{d}{q}+s<\sigma_++2.
\end{equation}
\end{lemma}

Recall that by solution we meant integral solution and therefore, we introduce a nonlinear mapping $\mathcal J$ given by
\[
   \mathcal{J}_{\varphi}[u](t):=e^{-t\mathcal{L}_a}\varphi+\int_0^te^{-(t-\tau)\mathcal{L}_a}F_{\alpha}(u(\tau))d\tau.
\]
A fixed point of this map would essentially be a solution to \eqref{heat equation}.
Next using Lemma \ref{nde}, 
 we establish the nonlinear estimates in Kato spaces with appropriate conditions on the parameters.

\begin{proposition}[Nonlinear estimate, sub-critical \& critical case]\label{nle}
Let  $\alpha>1$, $\gamma\in\R$ satisfy \eqref{E1} , \eqref{E2}.
Let $s\in \R$ and $q\in (1,\infty)$ satisfy 
\begin{equation}\label{ac0}
	\tau=s+\frac{d}{q}\le\tau_c= \frac{2+\gamma}{\alpha-1}.
\end{equation} 
Let $k,p$ satisfy
\begin{equation}\label{ac1}
\frac{\gamma}{\alpha-1}\le k,\qquad \alpha<p<\infty
\end{equation}
\begin{equation}\label{ac4}
\frac{s+\gamma}{\alpha}\le k
\end{equation}
\begin{equation}\label{ac2}
\sigma_-<k+\frac{d}{p}<\frac{\sigma_++2+\gamma}{\alpha}
\end{equation}
 \begin{equation}\label{ac3}
\frac{1}{\alpha}(\frac{d}{q}+s+\gamma)<\begin{cases}
\frac{d}{p}+k\le\tau \quad \text{ if } \ \tau<\tau_c\\
\frac{d}{p}+k<\tau \quad  \text{ if } \ \tau=\tau_c
\end{cases}.
\end{equation} 
Then for any $u,v\in \mathcal{K}^{p,q}_{k,s}(T)$ we have
\begin{equation*}
\begin{rcases}\|\mathcal{J}_{\varphi}[u]-\mathcal{J}_{\varphi}[v]\|_{\mathcal{K}^{p,q}_{k,s}(T)}\\
\|\mathcal{J}_{\varphi}[u]-\mathcal{J}_{\varphi}[v]\|_{\mathcal{K}^{q,q}_{s,s}(T)}
\end{rcases}
\lesssim T^{\frac{\alpha-1}{2}(\tau_c-\tau)}(\|u\|_{\mathcal{K}^{p,q}_{k,s}(T)}^{\alpha-1}+\|v\|_{\mathcal{K}^{p,q}_{k,s}(T)}^{\alpha-1})\|u-v\|_{\mathcal{K}^{p,q}_{k,s}(T)}.
\end{equation*}
\end{proposition}

\begin{remark}\label{rr}
Note that $\|v\|_{\mathcal{K}^{q,q}_{s,s}(T)}=\sup_{0\le t\le T}\|v(t)\|_{L_s^q}$.
\end{remark}


\begin{remark}\label{hy3}
First inequality in \eqref{ac2}, last inequality in
\eqref{ac3} and \eqref{ac0} imposes the condition $\sigma_-<\tau_c$.  This is equivalent with \begin{equation}\label{E2'}
\sigma_-<\frac{2+\gamma}{\alpha-1}\Longleftrightarrow\begin{cases}
\alpha<1+\frac{2+\gamma}{\sigma_-}\quad \text{if }\sigma_->0\\
0<2+\gamma\quad \text{if }\sigma_-=0\\
\alpha>1+\frac{2+\gamma}{\sigma_-}\quad \text{if }\sigma_-<0,
\end{cases}
\end{equation}  which is confirmed by \eqref{E1}, \eqref{E2} (using the fact $\sigma_->0\Leftrightarrow a<0$ and $\sigma_-<0 \Leftrightarrow a>0$ and $\sigma_-=0$ if $a=0$). 
\end{remark}

\begin{remark}
Note that \eqref{ac2} imposes the condition\begin{equation*}
\sigma_-<\frac{\sigma_++2+\gamma}{\alpha}\Longleftrightarrow\begin{cases}
\alpha<\frac{\sigma_+}{\sigma_-}+\frac{2+\gamma}{\sigma_-}\ \text{if }\sigma_->0\\
0<\sigma_++2+\gamma\ \text{if }\sigma_-=0\\
\alpha>\frac{\sigma_+}{\sigma_-}+\frac{2+\gamma}{\sigma_-}\ \text{if }\sigma_-<0
\end{cases}
\end{equation*}and this is implied by \eqref{E2'} as $\frac{\sigma_+}{\sigma_-}>1$ for $\sigma_->0$ and $\frac{\sigma_+}{\sigma_-}<0$ for $\sigma_-<0$.  
\end{remark}

 Before proving Proposition \ref{nle},  we  prove a technical lemma as an application of Theorem \ref{wle}.

\begin{lemma}\label{nde}Assume $\alpha\ge1$ and let $p\in(\alpha,\infty)$, $r\in(1,\infty)$, $l,k\in\R$ and 
\begin{equation}\label{f0}
\sigma_-<\frac{d}{r}+l,\quad \frac{d}{p}+k<\frac{\sigma_++2+\gamma}{\alpha},\qquad\gamma\le\alpha k- l+\min(\frac{\alpha d}{p}-\frac{d}{r},0).
\end{equation}
then for $t>0$ and $\varphi,\psi\in L_k^p(\rd)$ we have
\[
\|e^{-t\la}[|\cdot|^{\gamma}|\{\varphi|^{\alpha-1} \varphi)-|\psi|^{\alpha-1} \psi\}]\|_{L_l^r}\lesssim t^{-\frac{d
 }{2}(\frac{\alpha}{p}-\frac{1}{r})-\frac{\alpha k-l-\gamma}{2}}(\|\varphi\|_{L^{p}_k}^{\alpha-1}+\|\psi\|_{L^{p}_k}^{\alpha-1})\|\varphi-\psi\|_{L^{p}_k}.
\]

\end{lemma}\begin{proof} Note that \eqref{f0} is equivalent with  
 \[
\sigma_-<\frac{d}{r}+l\leq\frac{d}{p/\alpha}+\alpha k-\gamma<\sigma_++2,\qquad l\le\alpha k-\gamma.
\] 
By  Theorem \ref{wle},  with $s_2=l, s_1=\alpha k-\gamma , q_2=r, q_1=\frac{p}{\alpha}$ we obtain
\begin{eqnarray*}
&&\|e^{-t\la}[|\cdot|^{\gamma}|\{\varphi|^{\alpha-1} \varphi)-|\psi|^{\alpha-1} \psi\}]\|_{L_l^r}\\
&\lesssim&t^{-\frac{d
 }{2}(\frac{\alpha}{p}-\frac{1}{r})-\frac{\alpha k-\gamma-l}{2}}\||\cdot|^{\gamma}(|\varphi|^{\alpha-1} \varphi-|\psi|^{\alpha-1} \psi)\|_{L_{\alpha k-\gamma}^{\frac{p}{\alpha}}}\\
 &= &t^{-\frac{d
 }{2}(\frac{\alpha}{p}-\frac{1}{r})-\frac{\alpha k-l-\gamma}{2}}\||\cdot|^{\alpha k}(|\varphi|^{\alpha-1}+|\psi|^{\alpha-1})|\varphi-\psi|\|_{\frac{p}{\alpha}}
 \\&= &t^{-\frac{d
 }{2}(\frac{\alpha}{p}-\frac{1}{r})-\frac{\alpha k-l-\gamma}{2}}\|[(|\cdot|^{ k}|\varphi|)^{\alpha-1}+(|\cdot|^k|\psi|)^{\alpha-1}][|\cdot|^k(\varphi-\psi)]\|_{\frac{p}{\alpha}}.\end{eqnarray*}
 By using $\frac{\alpha}{p}=\frac{\alpha-1}{p}+\frac{1}{p}$ and Holder' inequality, the above quantity is dominated by
 \begin{eqnarray*} 
 & &t^{-\frac{d
 }{2}(\frac{\alpha}{p}-\frac{1}{r})-\frac{\alpha k-l-\gamma}{2}}\|(|\cdot|^k|\varphi|)^{\alpha-1}+(|\cdot|^k|\psi|)^{\alpha-1}\|_{L^{\frac{p}{\alpha-1}}}\||\cdot|^k(\varphi-\psi)\|_{L^p}\\
 &\lesssim &t^{-\frac{d
 }{2}(\frac{\alpha}{p}-\frac{1}{r})-\frac{\alpha k-l-\gamma}{2}}(\|\varphi\|_{L_k^p}^{\alpha-1}+\|\psi\|_{L_k^p}^{\alpha-1})\|\varphi-\psi\|_{L_k^p}
\end{eqnarray*}which completes the proof.
\end{proof}
\begin{proof}[\textbf{Proof of Proposition \ref{nle}}] Let us first establish two claims:\\
{\it Claim I:} Let $\beta=\beta(d,k,s,p,q)$ be as in Definition \ref{kato}. Then \begin{equation}\label{AC1}
	\beta\alpha<1
	\end{equation}
{\it Proof of Claim I:}	Note that $s+ \frac{d}{q}=\tau\le\tau_c= \frac{2+\gamma}{\alpha-1}$ implies $(s+\frac{d}{q})\alpha-2\le s+\frac{d}{q}+\gamma$. First inequality in  \eqref{ac3} says $s+\frac{d}{q}+\gamma<(\frac{d}{p}+k)\alpha$.  
	 Thus  
	 \[(s+\frac{d}{q})\alpha-2<(\frac{d}{p}+k)\alpha\Longleftrightarrow\left[\frac{d}{q}+s-\frac{d}{p}-k\right]\alpha<2\Longleftrightarrow\eqref{AC1}.
	 \]
	{\it Claim II:} 
	\begin{equation}\label{AC2}
	\frac{d}{2}(\frac{\alpha}{p}-\frac{1}{q})+\frac{\alpha k-\gamma-s}{2}<1
	\end{equation}
{\it Proof of Claim II:}		 For the subcritical case $\tau<\tau_c$ we have
	{\it Proof of claim:}\begin{eqnarray*}
\frac{d}{2}(\frac{\alpha}{p}-\frac{1}{q})+\frac{\alpha k-\gamma-s}{2}&\le &\frac{d}{2}(\frac{\alpha}{p}-\frac{1}{p})+\frac{\alpha k-\gamma-k}{2}\\
&=&\frac{1}{2}(\frac{d}{p}+k)(\alpha-1)-\frac{\gamma}{2}\\
&\le&\frac{1}{2}(\frac{d}{q}+s)(\alpha-1)-\frac{\gamma}{2}<1;
\end{eqnarray*}
 where in the first and third inequalities we used $\frac{d}{p}+k\le\tau=\frac{d}{q}+s$ and in the last step we used $\tau<\tau_c$.  Proof for the case $\tau=\tau_c$, we only need to make the first, third nonstrict inequalities by strict inequalities (using $\frac{d}{p}+k<\tau=\frac{d}{q}+s$) and last strict inequality by equality (using $\tau=\tau_c$). This proves Claim II.
 
	Now note that \eqref{ac1}, \eqref{ac2} implies \eqref{f0} for
	 $(p, r,l,s)=(p,p,k,k)$.  By Lemma \ref{nde} with $(p, r,l,s)=(p,p,k,k)$ and \eqref{Fcon} we have 
\begin{eqnarray}\label{h2}
   &&\|\mathcal{J}_{\varphi}[u]-\mathcal{J}_{\varphi}[v]\|_{L^{p}_k}\nonumber\\
   &\lesssim&\int_0^t\|e^{-(t-\tau)\mathcal{L}_a}[|x|^{\gamma}(|u|^{\alpha-1}u-|v|^{\alpha-1}v)(\tau)]\|_{L^{p}_k}d\tau\nonumber\\
   &\lesssim&\int_0^t(t-\tau)^{-\frac{d(\alpha-1)}{2p}-\frac{1}{2}\{(\alpha-1)k-\gamma\}}(\|u(\tau)\|_{L^{p}_k}^{\alpha-1}+\|v(\tau)\|_{L^{p}_k}^{\alpha-1})\|u(\tau)-v(\tau)\|_{L^{p}_k}d\tau\nonumber\\
   &
   \lesssim&(\|u\|_{\mathcal{K}^{p,q}_{k,s}(T)}^{\alpha-1}+\|v\|_{\mathcal{K}^{p,q}_{k,s}(T)}^{\alpha-1})\|u-v\|_{\mathcal{K}^{p,q}_{k,s}(T)}\int_0^t(t-\tau)^{-\frac{d(\alpha-1)}{2p}-\frac{1}{2}\{(\alpha-1)k-\gamma\}}\tau^{-\beta\alpha}d\tau,
\end{eqnarray} where the last inequality is due to the fact $u,u\in \mathcal{K}^{p,q}_{k,s}(T)$. 
Recall $\tau_c= \frac{2+\gamma}{\alpha-1}$ and $B(x,y):=\int_0^1\tau^{x-1}(1-\tau)^{y-1}d\tau$ is convergent if $x,y>0$.  Taking \eqref{ac0},  \eqref{ac3}, \eqref{AC1} into account,   
note that 
 the last time-integral in \eqref{h2} is bounded by
\begin{eqnarray*}
&&t^{1-\frac{d(\alpha-1)}{2p}-\frac{1}{2}\{(\alpha-1)k-\gamma\}-\alpha\beta}\int_0^1(1-\tau)^{-\frac{d(\alpha-1)}{2p}-\frac{1}{2}\{(\alpha-1)k-\gamma\}}\tau^{-\alpha\beta} d\tau\\
&=&t^{\frac{(\alpha-1)}{2}\left(\tau_c-\tau\right)}t^{-\beta}B\left(\frac{(\alpha-1)}{2}\big(\tau_c-\frac{d}{p}-k\big),1-\alpha\beta\right)< \infty.
\end{eqnarray*}
This together with \eqref{h2} implies the first part of the result.

 Note that \eqref{ac2}, \eqref{ac3} implies \eqref{f0} for
	 $(p, r,l,s)=(p,q,k,s).$  So by Lemma \ref{nde} with $(p, r,l,s)=(p,q,k,s),$ we have
\begin{eqnarray}\label{h3}
&&\|\mathcal{J}_{\varphi}[u](t)-\mathcal{J}_{\varphi}[v](t)\|_{L^{q}_s}\nonumber\\
&\lesssim&\int_0^t\|e^{-(t-\tau)\mathcal{L}_a}[|\cdot|^{-\gamma}(|u|^{\alpha-1}u-|v|^{\alpha-1}v)(\tau)]\|_{L^{q}_s} d\tau\nonumber\\
&\lesssim&\int_0^t(t-\tau)^{-\frac{d}{2}(\frac{\alpha}{p}-\frac{1}{q})-\frac{\alpha k-\gamma-s}{2}}(\|u(\tau)\|_{L^{p}_k}^{\alpha-1}+\|v(\tau)\|_{L^{p}_k}^{\alpha-1})\|u(\tau)-v(\tau)\|_{L^{p}_k}d\tau\nonumber\\
&\lesssim&(\|u\|_{\mathcal{K}^{p,q}_{k,s}(T)}^{\alpha-1}+\|v\|_{\mathcal{K}^{p,q}_{k,s}(T)}^{\alpha-1})\|u-v\|_{\mathcal{K}^{p,q}_{k,s}(T)}\int_0^t(t-\tau)^{-\frac{d}{2}(\frac{\alpha}{p}-\frac{1}{q})-\frac{\alpha k-\gamma-s}{2}}\tau^{-\alpha\beta}d\tau.
\end{eqnarray}
The last integral  is bounded by
\begin{eqnarray}\label{h4}
&&t^{1-\frac{d}{2}(\frac{\alpha}{p}-\frac{1}{q})-\frac{\alpha k-\gamma-s}{2}-\alpha\beta}\int_0^1(1-\tau)^{-\frac{d}{2}(\frac{\alpha}{p}-\frac{1}{q})-\frac{\alpha k-\gamma-s}{2}}\tau^{-\alpha\beta}d\tau\nonumber\\
&=&t^{\frac{(\alpha-1)}{2}(\tau_c-\tau)}\int_0^1(1-\tau)^{-\frac{d}{2}(\frac{\alpha}{p}-\frac{1}{q})-\frac{\alpha k-\gamma-s}{2}}\tau^{-\alpha\beta}d\tau,
\end{eqnarray}which is finite  in view of  \eqref{AC1} and \eqref{AC2}.
Now \eqref{h3} and \eqref{h4} implies the second part of the result.
\end{proof}

\begin{remark}[Hypotheses of Proposition \ref{nle}]\label{hpcp} \
\begin{itemize}
\item Condition \eqref{ac0} and last inequality in \eqref{ac3} are used to make sure the beta functions $B(x,y)$ is finite for various choices of $x,y$.
\item Conditions in \eqref{ac1}, \eqref{ac2} \eqref{ac3} are used to invoke Lemma \ref{nde} with $(p, r,l,s)=(p,p,k,k)$ and with $(p, r,l,s)=(p,q,k,s)$. 
\end{itemize}
\end{remark}

In the next result, we prove that there exists parameter $p,k$ such that \eqref{condition:p,k} in Lemma \ref{Kle} and \eqref{ac1}, \eqref{ac2}, \eqref{ac3} in Proposition \ref{nle} are satisfied.
\begin{lemma}\label{exi}
Assume \eqref{E1}, \eqref{E2}. Let  $\frac{\gamma}{\alpha-1}\le s$, $\sigma_--\frac{d}{\alpha}<s$ and $q\in (1,\infty)$ satisfy $\sigma_-<\frac{d}{q}+s<\sigma_++2$. 
Then there exist $k\in \R$ and $p\in (\alpha,\infty)$ satisfying
hypothesis \eqref{condition:p,k} of Lemma \ref{Kle},  and hypotheses \eqref{ac1}, \eqref{ac2}, \eqref{ac3}  of Proposition \ref{nle}. 
If we further assume $\tau<\tau_c$, \eqref{f01}, we can choose $p=q$ and $k=s$.
\end{lemma}
\begin{proof}
We need \begin{equation}\label{AA1}
	\sigma_{-}<\frac{d}{p}+k<\frac{d}{q}+s<-\gamma+(\frac{d}{p}+k)\alpha<\sigma_++2,
\end{equation} and 
\begin{equation}\label{AA2}
\frac{s+\gamma}{\alpha}\le k\le s.
\end{equation}

Now \eqref{AA1} follows if we chosse $\frac{d}{p}+k$ so that
\[
\max(\sigma_{-},\frac{\tau+\gamma}{\alpha})<\frac{d}{p}+k<\min(\frac{\sigma_++2+\gamma}{\alpha},\tau).
\]
Choose $k$ such that $$\max(\sigma_{-}-\frac{d}{\alpha},\frac{s+\gamma}{\alpha})<k<\min(\frac{\sigma_++2+\gamma}{\alpha},s)
$$so that \eqref{AA2} is satisfied. Then choose $p$ so that
\[
\max(\sigma_--k,\frac{\tau+\gamma}{\alpha}-k,0)<\frac{d	}{p}<\min(\frac{\sigma_++2+\gamma}{\alpha}-k,\frac{d}{q}+s-k,\frac{d}{\alpha})
\]which is possible as $\sigma_-<\frac{\sigma_++2+\gamma}{\alpha}$ as a consequence of \eqref{E2}.
This completes the proof.

The furthermore more part is clear.
\end{proof}

As we are done with linear estimate Lemma \ref{Kle}
and nonlinear estimate \ref{nle} and existence of parameter $p,k$ we are in a position to prove the following well-posedness result which implies Theorem \ref{lwp0}.
\begin{theorem}[Local well-posedness in the subcritical weighted Lebesgue space]\label{lwp2}
Let 
$\alpha>1$, $\gamma\in \R$ satisfy \eqref{E1} , \eqref{E2}.
Let $s\in \R$, $q\in (1,\infty)$ satisfy the subcriticality 
condition defined in \eqref{cri} and
\begin{equation}
\label{LWP.c.ql}
\frac{\gamma}{\alpha-1}\le s,\qquad \sigma_--\frac{d}{\alpha}<s.
\end{equation}
Let $k\in \R$ and $p\in (\alpha,\infty)$ satisfy 
hypothesis \eqref{condition:p,k} of Lemma \ref{Kle},  and hypotheses \eqref{ac1}, \eqref{ac2}, \eqref{ac3}  of Proposition \ref{nle}.
Then the Cauchy problem \eqref{heat equation} is locally well-posed in $L^q_s(\R^d)$ for arbitrary data $u_0\in L_s^q(\R^d)$. More precisely, the following assertions hold. 
\begin{enumerate}
\item\label{lwp2a} (Existence) 
	For any $u_0 \in L^{q}_s(\R^d),$ there exist a positive number 
	$T$ 
	 and an 
	$L^{q}_{s}(\R^d)$-integral solution $u$ to \eqref{heat equation} satisfying 
	\begin{equation}\nonumber	
		\|u\|_{\mathcal{K}^{p,q}_{k,s}(T)} 
			\le 2 \|e^{-t\mathcal{L}_a} u_0\|_{\mathcal{K}^{p,q}_{k,s}(T)}. 
	\end{equation}
	Moreover, the solution can be extended to the maximal interval 
	$[0,T_m)$. 
\item\label{lwp2b} (Uniqueness in ${\mathcal{K}}^{p,q}_{k,s}(T)$) 
	Let $T>0.$ If $u, v \in {\mathcal{K}}^{p,q}_{k,s}(T)$ satisfy 
	\eqref{integraleq} with $u(0) = v(0)=u_0$, then $u=v$ on $[0,T].$ 
\item (Continuous dependence on initial data)
	For any initial data $\varphi$ and $\psi$ in $L^{q}_{s}(\R^d),$  
	let $T(\varphi)$ and $T(\psi)$ be the corresponding existence time given by part  \eqref{lwp2a}. 
	Then there exists a constant $C$ depending on $\varphi$ and $\psi$ such that 
	the corresponding solutions $u$ and $v$ satisfy 
	\begin{equation}\nonumber	
		\|u-v\|_{L^\infty(0,T;L^{q}_{s}) \cap {\mathcal{K}}^{p,q}_{k,s}(T)} 
		\le C T^{\frac{\alpha-1}{2}(\tau_c-\tau)} 
						\|u_0-v_0\|_{L^{q}_{s}}
	\end{equation}
	for $T< \min\{T(u_0), T(v_0)\}.$  
\item\label{blc1} (Blow-up criterion in subcritical case $\tau<\tau_c$)) If $T_m<\infty,$ 
	then $\displaystyle \lim_{t\uparrow T_m}\|u(t)\|_{L^{q}_{s}}=\infty.$ 
	Moreover, the following lower bound of blow-up rate holds: 
	there exists a positive constant $C$ independent of $t$ such that 
	\begin{equation}\label{t:HH.LWP:Tm}
		\|u(t)\|_{L^{q}_{s}} 
		\gtrsim {(T_m - t)^{-\frac{\alpha-1}{2}(\tau_c-\tau)}} 
	\end{equation}
	for $t\in (0,T_m)$.
	\item\label{blc2} (Blow-up criterion in critical case $\tau=\tau_c$) 
If $u$ is an $L^q_{s}(\R^d)$-integral solution constructed in the assertion \eqref{lwp2a} and 
$T_m<\infty,$ then $\|u\|_{\mathcal{K}^{p,q}_{k,s}(T_m)}=\infty.$
	\item (Small data global existence 
	in critical case $\tau=\tau_c$)
There exists $\epsilon_0>0$ depending only on $d,\gamma,\alpha,q$ and $s$ such that if $u_0 \in \mathcal{S}'(\R^d)$ satisfies 
$\|e^{-t\mathcal{L}_a}u_0\|_{\mathcal{K}^{p,q}_{k,s}}<\epsilon_0$ (or $\|u_0\|_{L_s^q}<\epsilon_0$ in view of Lemma \ref{Kle}), then $T_m=\infty$ and 
$\|u\|_{\mathcal{K}^{p,q}_{k,s}} \le 2\epsilon_0$. 
\end{enumerate}
\end{theorem}


\begin{proof}[\textbf{Proof  of Theorem \ref{lwp2}}]
\underline{Existence in Kato space $\mathcal{K}_{k,s}^{p,q}(T)$:} Define $$B_M^T:=\{u\in\mathcal{K}_{k,s}^{p,q}(T)  : \|u\|_{\mathcal{K}_{k,s}^{p,q}(T)}\leq M\}$$with the metric
\[
d(u,v)=:  \|u-v\|_{\mathcal{K}_{k,s}^{p,q}(T)}.
\]
Then by Lemma \ref{Kle}, Proposition \ref{nle}, for $u,v\in B_M^T$ we have 
\begin{eqnarray}
\|\mathcal{J}_{u_0}[u]\|_{\mathcal{K}_{k,s}^{p,q}(T)}
&\le&\|e^{-t\la}u_0\|_{\mathcal{K}_{k,s}^{p,q}(T)}+cT^{\frac{\alpha-1}{2}(\tau_c-\tau)}M^{\alpha}\label{eqA}
\end{eqnarray}
and
\begin{eqnarray}
&&\|\mathcal{J}_{u_0}[u]-\mathcal{J}_{v_0}[v]\|_{\mathcal{K}_{k,s}^{p,q}(T)}\nonumber\\
&\le&\|e^{-t\la}(u_0-v_0)\|_{\mathcal{K}_{k,s}^{p,q}(T)}+cT^{\frac{\alpha-1}{2}(\tau_c-\tau)}M^{\alpha-1}\|u-v\|_{\mathcal{K}^{p,q}_{k,s}(T)}\label{eqB}
\end{eqnarray}
{\it   Subcritical case $\tau<\tau_c$:}
Using \eqref{eqA}, \eqref{eqB} and choosing $M=2\|e^{-t\la}u_0\|_{\mathcal{K}_{k,s}^{p,q}(T)}$ and $T>0$ small enough so that $cT^{\frac{\alpha-1}{2}(\tau_c-\tau)}M^{\alpha-1}\le\frac{1}{2}$, we find $\mathcal{J}_{u_0}$ is a contraction in $B_M^T$ sub-critical case (and hence we have existence of unique solution $u\in B_M^T$). This proves \eqref{lwp2a}, \eqref{lwp2b}.\\\\
{\it   Critical case $\tau=\tau_c$:} Note that using a density argument we have 
\begin{equation*}
\lim_{T\to0}\|e^{-t\la}u_0\|_{\mathcal{K}_{k,s}^{p,q}(T)}=0.
\end{equation*}Thus we choose $T>0$ so that $M:=2\|e^{-t\la}u_0\|_{\mathcal{K}_{k,s}^{p,q}(T)}$
and  $cM^{\alpha-1}<\frac{1}{2}$ where $c$ as in \eqref{eqA}.
Then by using \eqref{eqA}, \eqref{eqB}  for $u,v\in B_M^T$ we have 
\begin{eqnarray*}
\|\mathcal{J}_{u_0}[u]\|_{\mathcal{K}_{k,s}^{p,q}(T)}
&\le&\|e^{-t\la}u_0\|_{\mathcal{K}_{k,s}^{p,q}(T)}+cM^{\alpha}\label{eqA1}\le \frac{M}{2}+\frac{M}2= M
\end{eqnarray*}
and
\begin{eqnarray*}
\|\mathcal{J}_{u_0}[u]-\mathcal{J}_{v_0}[v]\|_{\mathcal{K}_{k,s}^{p,q}(T)}&\le&\|e^{-t\la}(u_0-v_0)\|_{\mathcal{K}_{k,s}^{p,q}(T)}+\frac12\|u-v\|_{\mathcal{K}^{p,q}_{k,s}(T)}\label{eqB1}
\end{eqnarray*} 
Thus $\mathcal{J}_{u_0}$ is a contraction in $B_M^T$. 
 This proves \eqref{lwp2a}. \\\\
\underline{Solution is in  $C([0,T),L_s^q(\rd))$:}
Using Lemma \ref{Kle}, Proposition \ref{nle} \begin{eqnarray*}
\|\mathcal{J}_{u_0}[u](t)\|_{L^{q}_s}
&\lesssim&\|u_0\|_{L^{q}_s}+M^{\alpha}T^{\frac{\alpha-1}{2}(\tau_c-\tau)}
\end{eqnarray*}
and
\begin{eqnarray*}\|\mathcal{J}_{u_0}[u](t)-\mathcal{J}_{v_0}[v](t)\|_{L^{q}_s}
&\lesssim&\|u_0-v_0\|_{L^{q}_s}+M^{\alpha-1}\|u-v\|_{\mathcal{K}^{p,q}_{k,s}(T)}T^{\frac{\alpha-1}{2}(\tau_c-\tau)}.
\end{eqnarray*}
Since $ \mathcal{J}_{u_0}[u]=u$, solution is indeed in $L^\infty([0,T);L_s^q(\rd))$ rest of the results  follows as in classical case.\\

Uniqueness, continuous dependency, blow-up, small data global existence are usual as in classical case.
\end{proof}

\begin{remark} \label{fr} The hypothesis  \eqref{f02}, \eqref{f01} are to make a integral functional map   contraction in Kato spaces.  
On the other hand conditions on $\alpha $ in \eqref{E2} are to make sure there exists $\tau$ satisfying  \eqref{f02}. 
\end{remark}

\begin{proof}[\textbf{Proof of Theorem \ref{lwp0}}]
The result follows from Theorem \ref{lwp2} and Lemma \ref{exi}. For the furthermore part, as in Lemma \ref{exi} we are choosing $p=q,k=s$, the uniqueness part follows from the uniqueness of fixed point in $B_T^M\subset\mathcal{K}^{q,q}_{s,s}(T) $ and Remark \ref{rr}.
\end{proof}
\begin{remark}\label{rR}
In \cite{chikami2022optimal}, Kato space $\mathcal{K}^{q,q}_{s,s}(T)$ was not used and hence the they did not achieve uniqueness in mere $C([0,T_m),L^q_s(\rd))$.
\end{remark}

\section{Finite time blow-up and nonexistence results}\label{S5}

In this section we establish that in the sub-critical and critical case there exists initial data for which solution established by Theorem \ref{lwp0} cannot be extended globally in time. Then blow-up alternative (see Theorem \ref{lwp2}) implies solution must blow-up in finite time. On the other hand for super-critical case, we shall prove that there exists data such that no local weak (hence integral) solution exists.

Before proving the above two we establish the following important lemma which will be used in both he proofs.

\begin{lemma}\label{ineq1}
Assume \eqref{scon0}. Let $u$ be a non-negative weak solution on $[0,T)$ to \eqref{heat equation} with initial data $u_0$. Let $\phi\in C^\infty_0(\mathbb R^d,[0,1])$ be such that $\phi=1$ on $B_{1/2}$ and supported in $B_1$. \
Then for $l\ge \max(3,\frac{2\alpha}{\alpha-1})$,  we have $$\int_{|x|<\sqrt{T}} u_0(x) \phi^l\left(\frac{x}{\sqrt{T}}\right)\, dx
 \lesssim T^{- \frac{2+\gamma}{2(\alpha-1)} + \frac{d}{2}}.$$
\end{lemma}
\begin{proof}
Let
\[
\psi_T(t,x)=\eta(\frac{t}{T})\phi(\frac{x}{\sqrt{T}}).
\]where $\eta \in C^\infty_0(\R,[0,1])$ is such that $\eta=1$ on $B_{1/2}$ and supported in $B_1$.
We note that for $l\ge3$ we have $\psi_T^l\in C^{1,2}([0,T)\times \R^d)$ and the estimate  
\begin{eqnarray}\label{n1}
|\partial_t \psi_T ^l(t,x)|+|\Delta\psi_T^l(t,x)|&\lesssim& T^{-1} \psi_T^{l-2}(t,x)\nonumber\\
&\lesssim& T^{-1} \psi_T^{\frac l \alpha}(t,x)
\end{eqnarray}by choosing \[l\ge\frac{2\alpha}{\alpha-1} \Longleftrightarrow\frac{l}{\alpha}\le l-2.\]  
We define a function $I:[0,T)\rightarrow \R_{\ge 0}$ given by
\[
  I(T):=\int_{[0,T)\times \{|x|<\sqrt{T}\}}|x|^{\gamma} u(t,x)^{\alpha} \, \psi_T^l(t,x) \, dtdx.
\]
We note that $I(T)<\infty$, since $u\in L_t^{\alpha}(0,T;L^{\alpha}_{\frac{\gamma}{\alpha},loc}(\R^d))$.
By using the weak form (\ref{weak}), non-negativity of $u$, the above estimate \eqref{n1}, H\"older's inequality and Young's inequality, the estimates hold:
\begin{eqnarray}\label{n0}
I(T) + \int_{|x|<\sqrt{T}} u_0(x) \phi^l\left(\frac{x}{\sqrt{T}}\right)\, dx
&=& \left|\int_{[0,T)\times \{|x|<\sqrt{T}\}}u(\partial_t \psi_T^l + \Delta \psi_T^l+a|x|^{-2}\psi_T^l )\,dt\,dx \right|\nonumber\\
& \le &\int_{[0,T)\times \{|x|<\sqrt{T}\}}(CT^{-1}+|a||x|^{-2})|u| \psi_T^{\frac{l}{\alpha}} dtdx\nonumber\\
&\le& CI(T)^{\frac{1}{\alpha}}K(T)^{\frac{1}{\alpha'}}\nonumber\\
&\le& \frac{1}{2}I(T)+CK(T),
\end{eqnarray}
where $1= \frac{1}{\alpha} + \frac{1}{\alpha'}$, i.e., $\alpha'=\frac{\alpha}{\alpha-1}$ and $K(T)$ is defined by
\[
K(T):=\int_{[0,T)\times\{|x|<\sqrt{T}\}}\{T^{-\alpha'}|x|^{-\frac{\gamma\alpha'}{\alpha}}+|a|^{\alpha'}|x|^{-\left(2+\frac{\gamma}{\alpha}\right)\alpha'}\}dxdt\sim T^{{-\frac{2+\gamma}{2(\alpha-1)}+\frac{d}{2}}}.
\]
The last equality holds only when \eqref{scon0} holds. 
Now from \eqref{n0}, we have
\begin{equation*}
\int_{|x|<\sqrt{T}} u_0(x) \phi^l\left(\frac{x}{\sqrt{T}}\right)\, dx\lesssim K(T)\sim T^{- \frac{2+\gamma}{2(\alpha-1)} + \frac{d}{2}}
\end{equation*}which completes the proof.
\end{proof}

\subsection{Finite time blow-up in critical and subcritical case}

The proof of this theorem is based on the arguments of \cite[Proposition 2.2, Theorem 2.3]{II-15} where lifespan of solution for nonlinear Schr\"odinger equation is studied.


\begin{proof}[{\bf Proof of Theorem \ref{blow-up large data}}]
Let $\lambda>0$ be a parameter. We take an initial data $u_0$ as $\lambda f$, where $f:\R^d\rightarrow \R_{\ge 0}$ is given by \begin{equation}
\label{initial data}
f(x) := 
\begin{cases}
	|x|^{-\beta} \quad & |x|\le 1,\\
	0 & \text{otherwise}
\end{cases}
\end{equation} with $\beta$ satisfying 
\begin{equation}\label{beta2'}
\beta< \min\left\{s+ \frac{d}{q},d\right\}.
\end{equation}
Then we see $u_0\in L^q_s (\R^d)$ and hence by Theorem \ref{lwp0}, 
we can define the maximal existence time $T_m=T_{m}(u_0)=T_m(\lambda f)$. Moreover the solution with initial data $\lambda f$ would be nonnegative as heat kernel is so.
Since $T_m(\lambda f)\le T_m^w(\lambda f)$, it follows from a change of variable and then  Lemma \ref{ineq1} that for any $0<T<T_m(\lambda f)$
\begin{eqnarray*}
 \lambda T^{\frac{d-\beta}{2}}\int_{|y|<1/\sqrt{T}}|y|^{-\beta}\phi^l(y)dy&=&\lambda \int_{|x|<1}|x|^{-\beta}\phi^l(\frac{x}{\sqrt{T}})dx\\
 &\le&\lambda \int_{\rd}|x|^{-\beta}\phi^l(\frac{x}{\sqrt{T}})dx\\
 &=&\lambda \int_{|x|<\sqrt{T}}|x|^{-\beta}\phi^l(\frac{x}{\sqrt{T}})dx\le C T^{- \frac{2+\gamma}{2(\alpha-1)} + \frac{d}{2}}
\end{eqnarray*}which implies
\begin{equation}\label{ineq2}
\lambda\le CL_T^{-1} T^{\frac{\beta}{2}- \frac{2+\gamma}{2(\alpha-1)} }
\end{equation}where $L_T=\int_{|y|<1/\sqrt{T}}|y|^{-\beta}\phi^l(y)dy$.

{\bf{Claim}}: There exists $\lambda_0$ such that if $\lambda>\lambda_0$, then $T_m(\lambda f)\le 4$.\\
Indeed, on the contrary we assume that $T_m(\lambda_j f)>4$ for a sequence $\lambda_j\to\infty$. Since $\beta<d$, we have $L_T<\infty$. The following estimates hold: 
\begin{eqnarray*}
  \lambda_j \le CL_4^{-1} 4^{\frac{\beta}{2}- \frac{2+\gamma}{2(\alpha-1)} }<\infty
\end{eqnarray*}
which a contradiction  and hence the claim is established.

Let $\lambda>\lambda_0$ and $0<T<T_m(\lambda f)\le4$ then again using \eqref{ineq2}
\[
\lambda\le CL_T^{-1}T^{ \frac{\beta}{2}- \frac{2+\gamma}{2(\alpha-1)} }\le CL_4^{-1}T^{ \frac{\beta}{2}- \frac{2+\gamma}{2(\alpha-1)} }
\]as $L_T$ is decreasing in $T$. By \eqref{beta2'} and the fact $\tau\le\tau_c$ we have $\kappa:= \frac{2+\gamma}{2(\alpha-1)}-\frac{\beta}{2}>0$ and so for all $T\in(0,T_m(\lambda f)$
\[
T\le c\lambda^{-\frac{1}{\kappa}}
\]which implies $T_m(\lambda f)\le c\lambda^{-\frac{1}{\kappa}}$.

Then the result follows from blowup criterion in Theorem \ref{lwp2}\eqref{blc1}. First point in Remark \ref{rfbl} follows from Theorem \ref{lwp2}\eqref{blc2}.
\end{proof}

\subsection{Nonexistence of weak solution in the supercritical case}

In this subsection we give a proof of Theorem \ref{t:nonex}. 
We only give a sketch of the proof. 
For the details, we refer to \cite[Proposition 2.4, Theorem 2.5]{II-15} where nonlinear Schr\"odinger equation is studied. 


\begin{proof}[{\bf Proof of Theorem \ref{t:nonex}}]
Let $T\in (0,1)$. Suppose that the conclusion of Theorem \ref{t:nonex} does not hold. 
Then there exists a positive weak solution $u$ on $[0,T)$ to (\ref{heat equation}) 
(See Definition \ref{weakdef}) with any initial data $u_0$ in particular for $f$ given by \eqref{initial data}
with $\beta$ satisfying
\begin{equation}\label{beta1}
	\frac{2+\gamma}{\alpha-1}<\beta< \min\left\{s+ \frac{d}{q},d\right\}.
\end{equation}
Note that such choice is possible as $\tau>\tau_c$ and \eqref{scon} i.e. $\alpha>\alpha_F(d,\gamma)$.
Now \eqref{beta1} implies $u_0 \in L^q_s (\R^d)\cap L^1_{loc}(\R^d)$. For $T<1$  we have using Lemma \ref{ineq1}
\begin{equation}\label{ineq2'}
\begin{aligned}
	\int_{|x|<\sqrt{T}} u_0(x) \phi^l\left(\frac{x}{\sqrt{T}}\right)\, dx
	& = T^{-\frac{\beta-d}{2}}	\int_{|y|<1} |y|^{-\beta} \phi^l(y)\, dx
	 = C T^{-\frac{\beta-d}{2}}.
\end{aligned}
\end{equation}
Combining Lemma \ref{ineq1} and \eqref{ineq2'}, we obtain
\begin{equation*}\label{contradiction}
0< C \le T^{\frac{\beta}{2} - \frac{2+\gamma}{2(\alpha-1)}} \to 0 \quad \text{as }T\to 0
\end{equation*}which leads to a contradiction, as 
$\beta$ satisfies 
\begin{equation*}
\frac{\beta}{2} - \frac{2+\gamma}{2(\alpha-1)} >0\quad \text{i.e.}\quad \beta > \frac{2+\gamma}{\alpha-1}.
\end{equation*}
  This completes the proof.
\end{proof}




\noindent
{\textbf{Acknowledgement}:} S Haque is thankful to DST--INSPIRE (DST/INSPIRE/04/2022/001457) \& USIEF--Fulbright-Nehru fellowship 
 for financial support. S Haque is also thankful to Harish-Chandra Research Institute \&  University of California, Los Angeles  for their excellent research facilities.  
 
\bibliographystyle{siam}
\bibliography{ref}
\end{document}